\documentclass{article}
\usepackage[utf8]{inputenc}
\usepackage{amsmath}
\usepackage{amssymb}
\usepackage{amsthm}
\usepackage{colonequals}
\usepackage{xcolor}
\usepackage{url}

\usepackage{tikz}
\usetikzlibrary{arrows, arrows.meta,  patterns, calc}
\usetikzlibrary{decorations.markings, intersections}

\newcommand{\hes}{\mathrm{hes}}
\newcommand{\ohes}{\vec{\mathrm{hes}}}
\newcommand{\opp}{\mathrm{opp}}
\newcommand{\lft}{\mathrm{left}}
\newcommand{\tgt}{\mathrm{tgt}}
\newcommand{\ew}{\mathrm{ew}}
\DeclareMathOperator{\pll}{\mathrm{polylog}}
\newcommand{\supp}{\mathrm{supp}}
\newcommand{\trans}[1]{\tau_{#1}}
\newcommand{\transplus}[1]{\tau^+_{#1}}

\newcommand{\PP}{\mathcal{P}}
\newcommand{\wtococ}[1]{\circ#1}

\tikzset{
vtx/.style={inner sep=1.2pt, outer sep=0pt, circle, fill=black,draw=black}, 
}

\tikzset{
vtxr/.style={inner sep=1.2pt, outer sep=0pt, circle, fill=red,draw=black}, 
vtxw/.style={inner sep=1.2pt, outer sep=0pt, circle, fill=white,draw=black},
midarrow/.style={decoration={
    markings,
    mark=at position 0.5 with {\arrow{Latex[length=1.5mm]}}},postaction={decorate}},
bigmidarrow/.style={decoration={
    markings,
    mark=at position 0.5 with {\arrow{Latex[length=2.5mm]}}},postaction={decorate}},    
rededge/.style={color=red, midarrow}, 
blackedge/.style={color=black}, 
}

\IfFileExists{fig-tikz.pdf}{}{%
\immediate\write18{pdflatex fig-tikz}
}
\IfFileExists{fig-tikz3colorable.pdf}{}{%
\immediate\write18{pdflatex fig-tikz3colorable}
}

\newtheorem{theorem}{Theorem}
\newtheorem{observation}[theorem]{Observation}
\newtheorem{lemma}[theorem]{Lemma}
\newtheorem{corollary}[theorem]{Corollary}

\newenvironment{proofout}{%
  \begin{proof}[Proof outline]%
}{%
  \end{proof}%
}

\title{Embedded graph 3-coloring and flows}
\author{
Caroline Bang\thanks{Department of Mathematics, Iowa State University, Ames, IA. E-mail: \url{cbang@iastate.edu}. Research of this author is supported in part by NSF grants DMS-1839918 and DMS-2152490.}
\and
Zden\v{e}k Dvo\v{r}\'ak\thanks{Computer Science Institute, Charles University, Prague, Czech Republic
E-mail: \url{rakdver@iuuk.mff.cuni.cz}. Supported by project 22-17398S (Flows and cycles in graphs on surfaces) of Czech Science Foundation.}
\and
Emily Heath\thanks{Department of Mathematics, Iowa State University, Ames, IA. E-mail: \url{eheath@iastate.edu}. Research of this author is supported in part by NSF grant DMS-1839918.}
\and
Bernard Lidick\'y\thanks{Department of Mathematics, Iowa State University, Ames, IA. E-mail: \url{lidicky@iastate.edu}. Research of this author is supported in part by NSF grant  DMS-2152490 and Scott Hanna fellowship.}
}
\date{\today}

\begin{document}

\maketitle

\begin{abstract}
A graph drawn in a surface is a \emph{near-quadrangulation} if the sum of the lengths of the faces different from $4$-faces
is bounded by a fixed constant.  We leverage duality between colorings and flows to design an efficient algorithm
for 3-precoloring-extension in near-quadrangulations of orientable surfaces.  Furthermore, we use this duality to
strengthen previously known sufficient conditions for 3-colorability of triangle-free graphs drawn in orientable surfaces.
\end{abstract}

\section{Introduction}

In general, it is NP-hard to decide whether a planar graph is 3-colorable~\cite{garey1979computers};
however, a well-known theorem of Gr\"otzsch~\cite{grotzsch1959} states
that every planar \emph{triangle-free} graph is 3-colorable.
This result motivated further exploration into sufficient conditions for 3-colorability of planar graphs,
see e.g.~\cite{borsurvey}, as well as for more general graph classes, such as graphs drawn on other surfaces.

A graph is \emph{$(k+1)$-critical} if it is not $k$-colorable, but all its proper subgraphs are $k$-colorable;
hence, $(k+1)$-critical graphs are exactly the minimal forbidden subgraphs for $k$-colorability.
Thus, Gr\"otzsch's theorem is equivalent to the fact that there are no planar triangle-free 4-critical graphs.
Gimbel and Thomassen~\cite{gimbel} extended this result by showing that a triangle-free graph drawn in the projective
plane is 4-critical if and only if it is a non-bipartite quadrangulation without separating 4-cycles.

It turns out that \emph{near-quadrangulations} play an important role in 3-colorability of triangle-free graphs
in any fixed surface.  For reasons that will become clear later (see Observation~\ref{obs-relev}), it is convenient to use the following definition.
For an integer $n$, let $q(n)$ be the number of integers $i$ such that $3|i$, $i\equiv n\pmod 2$, and $|i|\le n$,
and let $b(n)$ be the largest such integer $i$.  Let $H$ be a graph with a 2-cell drawing in a surface
and let $F(H)$ denote the set of faces of $H$.  For a face $f\in F(H)$, let $|f|$ denote the length of the closed walk in $H$ that bounds $f$.
We let
\begin{align*}
q^\star(H)&=\prod_{f\in F(H)} q(|f|)\text{ and}\\
b^\star(H)&=1+\sum_{f\in F(H)} b(|f|).
\end{align*}
Note that $q(4)=1$ and $b(4)=0$, and thus if $H$ is a quadrangulation, then $q^\star(H)=b^\star(H)=1$.
We say that $H$ is an \emph{$a$-near-quadrangulation} if $b^\star(H)\le a$.  Dvo\v{r}\'ak, Kr\'al' and Thomas~\cite{trfree4} proved the following key result.
\begin{theorem}[Dvo\v{r}\'ak, Kr\'al' and Thomas~\cite{trfree4}]\label{thm-trfree4}
For every surface $\Sigma$ of Euler genus $g$, there exists a positive integer $a_\Sigma=O(g)$ such that every
4-critical triangle-free graph $H$ drawn in $\Sigma$ satisfies at least one of the following conditions:
\begin{itemize}
\item The drawing of $H$ is not 2-cell, or
\item $H$ contains a non-contractible 4-cycle, or
\item $H$ is an $a_\Sigma$-near-quadrangulation.
\end{itemize}
\end{theorem}
Dvo\v{r}\'ak, Kr\'al' and Thomas~\cite{trfree6} also gave a linear-time algorithm to 3-color near-quadrangulations
and combined these results in a linear-time algorithm to decide 3-colorability of triangle-free graphs drawn in any
fixed surface~\cite{trfree7}.

The algorithm of~\cite{trfree6} for 3-coloring near-quadrangulations is quite complicated;
it uses cutting and precoloring arguments to transform the input instance into a generic one (where 
there are no short non-contractible cycles) at the cost of introducing precolored vertices incident with pairwise
distant faces and then characterizes precoloring extension in such a generic instance by a topological criterion.
The downside of this approach is that it is quite non-explicit, and
the complexity and large multiplicative constants make it unusable in practice.  Hence, it is interesting to investigate
alternative approaches.

To establish one of the basic cases (plane graph with vertices incident with the outer faces precolored and with all other faces
of length four),  Dvo\v{r}\'ak, Kr\'al' and Thomas~\cite{trfree6} used the duality with nowhere-zero flows.
Dvo\v{r}\'ak and Lidick\'y~\cite[Lemma~4 and the remarks after it]{dl2015} explored this connection in more detail
and gave the following algorithm.
\begin{theorem}[Dvo\v{r}\'ak and Lidick\'y~\cite{dl2015}]
There exists an algorithm that, given a simple $n$-vertex plane graph $H$, decides whether $H$ is 3-colorable (and finds a 3-coloring if it exists)
in time $O(q^\star(H)b^\star(H)n)$.  Moreover, this algorithm can also decide whether a precoloring of the vertices incident
with the outer face of $H$ extends to a 3-coloring of $H$.
\end{theorem}
The algorithm is based on at most $q^\star(H)$ invocations of a maximum flow algorithm, and thus it is easy to implement in practice.
In a similar vein, Dvo\v{r}\'ak and Pek\'arek~\cite{dpek} considered plane graphs with two precolored faces, showing the following
result.
\begin{theorem}[Dvo\v{r}\'ak and Pek\'arek~\cite{dpek}]\label{thm-cyl}
There exists an algorithm that, given a simple $n$-vertex plane graph $H$ and a precoloring $\psi$ of the vertices incident with two faces of $H$, decides
whether $\psi$ extends to a $3$-coloring of $H$ (and finds such a 3-coloring if it exists)
in time $O(q^\star(H)b^\star(H)n)$.
\end{theorem}
Dvo\v{r}\'ak and Pek\'arek~\cite{dpek} used their result to obtain a practical algorithm for 3-coloring near-quadrangulations of the
torus with bounded edgewidth, as well as for deciding 3-colorability of triangle-free toroidal graphs.

Our main result is a far-reaching generalization of this approach:
\begin{itemize}
\item We consider graphs drawn in any orientable surface and do not put any restriction on the edgewidth.
\item We allow an arbitrary subset of the vertices to be precolored.
\end{itemize}
\begin{theorem}\label{thm-main}
Let $\Sigma$ be an orientable surface of Euler genus $g$.  There exists a function $\gamma$ and an algorithm that,
given a simple $n$-vertex graph $H$ with a 2-cell drawing in $\Sigma$ and a precoloring $\psi$ of a subset $S$ of its vertices, decides
whether $\psi$ extends to a $3$-coloring of $H$ (and finds such a 3-coloring if it exists)
in time
$$O\left(q^\star(H)b^\star(H)n+q^\star(H)\cdot\min\genfrac{(}{)}{0pt}{}{n^g(n^2+|S|^3),}{\gamma(|S|)n^2\pll n}\right).$$
\end{theorem}
It is remarkable that in Theorem~\ref{thm-main} we get a polynomial algorithm even if an arbitrary number of vertices is precolored.
In contrast, deciding 3-colorability of bipartite graphs with only three precolored vertices is
NP-complete, as shown by Kratochv\'il and Seb\H{o}~\cite{kratochvil1997coloring}, and deciding 3-colorability in planar triangle-free
graphs of maximum degree four with (unbounded number of) precolored vertices was shown to be NP-complete by Monnot~\cite{monnot2006note}.

Compared to the algorithm of Dvo\v{r}\'ak, Kr\'al' and Thomas~\cite{trfree6}, our algorithm
is much more practical (it is based on a simple combination of off-the-shelf algorithms---maximum flow, shortest path,
integer programming in bounded dimension) and allows one to precolor an arbitrary number of vertices.
On the other hand, the time complexity of Dvo\v{r}\'ak et al.~\cite{trfree6} algorithm is linear in the number of vertices
of the input graph.

Compared with Theorem~\ref{thm-cyl}, we offer a worse dependence on the number of precolored vertices.
However, the algorithm from Theorem~\ref{thm-main} can be adjusted so that its complexity
does not depend on $|S|$, but only on the number of components of $H[S]$, thus bridging this gap.
This is achieved by contracting the edges between precolored vertices and adjusting the flow constraints
according to the flow amount forced on the duals of these edges by the precoloring; see~\cite[Observation 17]{dpek} for
a precise explanation of the idea, which can be easily adapted to our setting.

A \emph{homomorphism} from a graph $H$ to a graph $C$ is a function $f:V(H)\to V(C)$ such that for every $uv\in E(H)$,
we have $f(u)f(v)\in E(C)$.  A 3-coloring of $H$ is equivalent to a homomorphism to $C_3$.
Theorem~\ref{thm-main} can be generalized to homomorphisms to odd cycles.  A motivation to study such homomorphisms arises
from their relation to the \emph{circular chromatic number}: The circular chromatic number $\chi_c(H)$ of a graph $H$ is
the minimum length of a circle for which there exists a mapping from $V(H)$ to open arcs of length $1$ in the circle
such that the arcs of adjacent vertices are disjoint.  Circular chromatic number is a refinement of the ordinary chromatic
number~\cite{vince}, in the sense that $\chi(H)=\lceil \chi_c(H)\rceil$ for every graph $H$.  And, for any positive integer $k$,
a graph has circular chromatic number at most $2+1/k$ if and only if it has a homomorphism to $C_{2k+1}$.

For a cycle $C$ of length $m$ and an integer $n$, let $q_C(n)$ be the number of integers $i$ such that $m|i$, $i\equiv n\pmod 2$, and $|i|\le n$,
and let $b_C(n)$ be the largest such integer $i$.
For a graph $H$ with a 2-cell drawing in a surface, let $q_C^\star(H)=\prod_{f\in F(H)} q_C(|f|)$
and $b_C^\star(H)=1+\sum_{f\in F(H)} b_C(|f|)$.
\begin{theorem}\label{thm-maingen}
Let $\Sigma$ be an orientable surface of Euler genus $g$ and let $C$ be an odd cycle.  There exists a function $\gamma$ and an algorithm that,
given a simple $n$-vertex graph $H$ with a 2-cell drawing in $\Sigma$ and a function $\psi:S\to V(C)$ from a subset $S$ of its vertices, decides
whether $\psi$ extends to a homomorphism from $H$ to $C$ (and finds such a homomorphism if it exists)
in time
$$O\left(q_C^\star(H)b_C^\star(H)n+q_C^\star(H)\cdot\min\genfrac{(}{)}{0pt}{}{n^g(n^2+|S|^3),}{\gamma(|S|)n^2\pll n}\right).$$
\end{theorem}

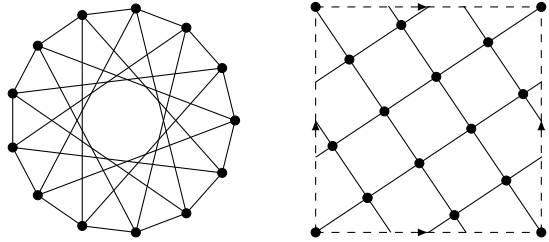
\begin{figure}
\begin{center}
\begin{tikzpicture}[scale=1.5]
\draw
\foreach \i in {0,1,...,17}{
(360/13*\i:1) coordinate (\i)
};
\foreach \i in {0,1,...,12}{
\pgfmathtruncatemacro{\ii}{\i+1}
\pgfmathtruncatemacro{\iii}{\i+5}
\draw
(\ii)--(\i)--(\iii)
(\i) node[vtx]{}
;
}
\end{tikzpicture}
\hskip 2em
\tikzset{
toruscut/.style={decoration={
    markings,
    mark=at position 0.5 with {\arrow{latex}}},postaction={decorate},dashed}, 
}
\begin{tikzpicture}
\draw[toruscut](0,0) -- (3,0);
\draw[toruscut](3,0) -- (3,3);
\draw[toruscut](0,0) -- (0,3);
\draw[toruscut](0,3) -- (3,3);
\begin{scope}
\clip(0,0) rectangle(3,3);
\draw[name path=l1](1,0) -- +(124:4);
\draw[name path=l2](2,0) -- +(124:4);
\draw[name path=l3](3,0) -- +(124:4);
\draw[name path=l4](4,0) -- +(124:4);
\draw[name path=a1](0,2) -- +(33.7:4);
\draw[name path=a2](0,1) -- +(33.7:4);
\draw[name path=a3](0,0) -- +(33.7:4);
\draw[name path=a4](1.5,0) -- +(33.7:3);
\draw (0,0) coordinate (intersection-1);
\foreach \l in {l1,l2,l3,l4}{
\foreach \a in {a1,a2,a3,a4}{
\draw[name intersections={of={\l} and {\a}}]
    (intersection-1) node[vtx]{ };
    }
}
\end{scope}
\draw (3,0) node[vtx]{} (3,3) node[vtx]{} (0,0) node[vtx]{} (0,3) node[vtx]{} ;
\end{tikzpicture}
\caption{Cayley graph $C(Z_{13}; 1, 5)$ and its drawing as a quadrangulation $Q_{13}$ of the torus.}\label{fig:Q13}
\end{center}
\end{figure}

Going back to the 3-coloring case, Hutchingson~\cite{locplanq} proved that for each orientable surface,
all graphs drawn in this surface with no odd-length faces and with sufficiently large edgewidth are 3-colorable
(the \emph{edgewidth} $\ew(H)$ of a graph $H$ drawn in a surface other than the sphere is the length of the shortest non-contractible cycle in $H$).
Dvo\v{r}\'ak, Kr\'al' and Thomas~\cite{trfree6} extended this claim to all triangle-free graphs.
Our argument gives a more general result with an explicit bound on edgewidth.  To state it, we need a few definitions.
Let $H$ be a graph drawn in an orientable surface $\Sigma$ and let $U$ be a subgraph of $H$.
\begin{itemize}
\item We say $U$ is \emph{flat} if $U$ does not contain any non-contractible cycle; or equivalently, there exists a disk $\Delta\subseteq \Sigma$ containing $U$.
We call the unique face of $U$ not contained in $\Delta$ the \emph{outer face} and all other faces \emph{internal faces}.
\item We say that $U$ \emph{captures non-4-faces of $H$} if every face of $H$ of length other than four is also
a face of $U$.
\item A graph $U'$ drawn in the plane is a \emph{planar quadrangulation extension} of $U$ if $U'$ is obtained from $U$ by quadrangulating its outer face.
More precisely, $U'$ is a planar quadrangulation extension of $U$ if there exists a homeomorphism $\theta$ from $\Delta$ to
a disk $\Delta'$ in the plane such that
\begin{itemize}
\item $\theta$ maps $U$ to a subgraph of $U'$,
\item $\theta$ maps internal faces of $U$ to faces of $U'$ contained in $\Delta'$, and
\item every face of $U'$ that is not the image of a face of $U$ under $\theta$ has length four.
\end{itemize}
\end{itemize}
A graph $H$ drawn in an orientable surface is \emph{locally $3$-colorable} if
every flat subgraph capturing non-4-faces of $H$ has a $3$-colorable planar quadrangulation extension.
\begin{theorem}\label{thm-ew}
Let $H$ be a simple graph with a 2-cell drawing in an orientable surface of Euler genus $g>0$ with
edgewidth at least $b^\star(H)+\Omega(g^{4/3})$.  Then $H$ is $3$-colorable if and only if it is locally 3-colorable.
\end{theorem}
If $H$ is a triangle-free graph, then every flat subgraph of $H$
has a triangle-free planar quadrangulation extension, which is 3-colorable by Gr\"otzsch's theorem~\cite{grotzsch1959}.
Hence, the condition of being locally 3-colorable is automatically satisfied.
In conjunction with Theorem~\ref{thm-trfree4}, this gives the following corollary.
\begin{corollary}\label{cor-3fc}
Let $H$ be a triangle-free graph drawn in an orientable surface of Euler genus $g>0$.
If $H$ has edgewidth $\Omega(g^{4/3})$, then $G$ is $3$-colorable.
\end{corollary}
\begin{proof}
Suppose for a contradiction that $H$ is not $3$-colorable.  Without loss of generality, we can assume that $H$ is 4-critical.
Moreover, we can assume that the drawing of $H$ is $2$-cell, as otherwise we can cut the surface along a non-contractible simple
closed curve contained in one of the faces and cap the resulting holes by disks, obtaining a drawing of $H$ in an orientable
surface of smaller genus, without decreasing the edgewidth.  We can also assume that the edgewidth is at least five.
By Theorem~\ref{thm-trfree4}, this implies that $b^\star(H)=O(g)$, and Theorem~\ref{thm-ew} shows that $H$ is 3-colorable.
\end{proof}
Let us remark that the existence of a lower bound on edgewidth guaranteeing 3-colorability of a triangle-free graph in an orientable surface
(together with an analogous, slightly more complicated result for non-orientable surfaces) has been proven in~\cite{trfree6},
without quantifying the dependence on the genus. For graphs on the torus, a more detailed analysis gives the following explicit bound.
\begin{corollary}\label{cor-torus}
Let $H$ be a simple graph with a 2-cell drawing on the torus.  If $H$ is triangle-free and has edgewidth at least $5+b^\star(H)$, then $H$ is $3$-colorable.
\end{corollary}
For quadrangulations of the torus, we have $b^\star(H)=1$, and Corollary~\ref{cor-torus} states that if the edgewidth is at least $6$,
then $H$ is 3-colorable.  This bound cannot be improved, as Archdeacon et al.~\cite{archdeacon2001chromatic} found a non-3-colorable quadrangulation $Q_{13}$ of
the torus with edgewidth five, see Figure~\ref{fig:Q13}.  Let us remark that Král' and Thomas~\cite{thomas2008coloring} proved that every non-3-colorable graph drawn of the torus
without odd faces contains $Q_{13}$ as a subgraph.

More generally, the results of Dvořák and Pekárek~\cite{dvopek} imply that if $H$ is a 4-critical triangle-free graph drawn on the torus, then
$b^\star(H)\le 13$; consequently, Corollary~\ref{cor-torus} implies that every triangle-free graph drawn on the torus with edgewidth at least $18$
is 3-colorable.  Let us remark that an exact characterization of 3-colorability of the triangle-free graphs drawn on the torus was obtained using computer-assisted enumeration~\cite{dpekchar}
and implies that every triangle-free graph drawn on the torus with edgewidth at least $6$ is 3-colorable.

Finally, for graphs without odd faces (and in particular for quadrangulations), we can strengthen the bound from Corollary~\ref{cor-3fc}.
\begin{theorem}\label{thm-genq}
For any orientable surface $\Sigma$ of Euler genus $g>0$, every graph $H$ with a $2$-cell drawing in $\Sigma$ of edgewidth $\Omega(g\log g)$
such that all faces have even length is $3$-colorable.
\end{theorem}

The rest of the paper is organized as follows:
\begin{itemize}
\item In Section~\ref{sec-prelim}, we introduce the basic notions from the homology theory and flow-coloring duality
and reduce the coloring of near-quadrangulations to the problem of finding circulations with prescribed homology.
In Section~\ref{sec-homol} we give an algorithm for this problem.
\item In Section~\ref{sec-allow}, we show that realizable homologies form a polytope and study its properties.
\item In Section~\ref{sec-nonem}, we argue that finding a suitable realizable homology reduces to finding an integer point in a related polytope
and give an algorithm for a special case arising when dealing with precolored vertices.  Combining the results
obtained till this point, we give the algorithm proving Theorem~\ref{thm-maingen}.
\item In Section~\ref{sec-hollow}, we state bounds on the width of polytopes with no integer pointsand show how edgewidth of the graph lower bounds the width of the polytopes relevant for its coloring.
Using these bounds, we prove Theorem~\ref{thm-genq}.
\item In Section~\ref{sec-loc3col}, we relate local 3-colorability to existence of nowhere-zero flows with boundary divisible by $3$
and prove Theorem~\ref{thm-ew} and Corollary~\ref{cor-torus}.
\end{itemize}
We finish in Section~\ref{sec-final} by some concluding remarks.

\section{Preliminaries}\label{sec-prelim}

We aim to find a coloring of a graph $H$ drawn in a surface by utilizing flows in its dual graph $G$
(the vertices of $G$ are the faces of $H$, and each edge of $H$ contributes an edge to $G$ joining the
two incident faces).  Let us remark that even though $H$ can be assumed without loss of generality to be a simple graph, $G$ may
have parallel edges and loops.  We will mostly focus on the properties of $G$, and when needed, we will use $G^\star$ to refer to $H$.

\begin{figure}
    \centering
    \includegraphics[page=2]{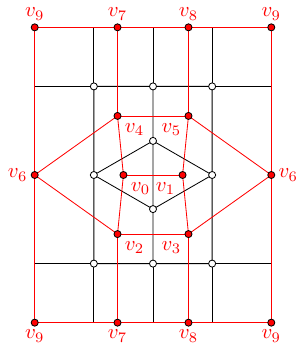}
    \caption{A graph $G$ drawn in red drawn on the torus (obtained by gluing the top edge of the picture with the bottom one, and the left edge with the right one). The dual $G^\star$ is drawn in black.}
    \label{fig:GGstar}
\end{figure}

We provide a number of figures to illustrate the concepts and results.  All figures should be interpreted as depicting a graph
drawn on the torus obtained by identifying the opposite sides of the rectangle.  We drawn the graph $G$ in red and the graph $G^\star$
in black.  For example, Figure~\ref{fig:GGstar} depicts a graph $G$ with $10$ vertices and with a double edge between vertices $v_6$ and $v_9$,
and its dual graph $G^\star$ (which is simple and also has $10$ vertices).  When we are asked to select faces in $G$, we will often
depict the corresponding vertices in $G^\star$ instead to avoid cluttering the picture.

\subsection*{Graphs on surfaces}

Let $G$ be a connected graph.  A \emph{drawing} $\eta$ of $G$ in a surface $\Sigma$ maps vertices of $G$ to pairwise distinct
points of $\Sigma$, each non-loop edge $e=uv\in E(G)$ to a simple curve in $\Sigma$ with ends $\eta(u)$ and $\eta(v)$,
and each loop $e=vv\in E(G)$ to a non-trivial simple closed curve containing $\eta(v)$, such that for each $e,e'\in E(G)$,
the intersection of the curves $\eta(e)$ and $\eta(e')$ consists only of the points $\eta(v)$ for vertices $v\in V(G)$ incident
with both $e$ and $e'$.  The \emph{faces} of the drawing are the maximal connected subsets of $\Sigma\setminus \bigcup_{e\in E(G)} \eta(e)$;
let $F(G)$ denote the set of faces of $G$.
The drawing is \emph{$2$-cell} if each face is homeomorphic to an open disk.  We only consider 2-cell drawings, since if
a connected graph has a non-2-cell drawing in a surface $\Sigma$, it also has a drawing in a surface of smaller genus.
In particular, any drawing of a connected graph in the sphere is $2$-cell.  We also restrict ourselves
to orientable surfaces (we briefly discuss non-orientable surfaces in Section~\ref{sec-final}).

\begin{figure}
\begin{center}
    \begin{tikzpicture}[thick]
    \fill[pattern=north west lines, pattern color=blue!30] (-.8, 1.3) -- (-.8, .4) -- (0,0) -- (3,0) -- (3.8, .4) -- (3.8, 1.3) -- (-.8, 1.3);
    \fill [black] (0,0) circle (.1) node [below] at (0,0) {$u$};
    \fill [black] (3,0) circle (.1) node [below] at (2.5,0) {$v = \tgt(h)$};
    \draw[red, <-] (0.1,0)-- (1.5,0) node [above] at (0.75, 0) {$\opp(h)$};
    \draw[blue, ->] (1.5,0)-- (2.9,0) node [above] at (2.25, 0) {$h$};
    \draw (0,0) -- (-.8, .4) -- (-.8, 1.3);
    \draw (3,0) -- (3.8, .4) -- (3.8, 1.3);
    \draw (0,0) -- (-.8, -.4) -- (-.8, -1.3);
    \draw (3,0) -- (3.8, -.4) -- (3.8, -1.3);
    \draw node at (1.5, 0.9) {$\lft(h)$};
    \end{tikzpicture}
    \end{center}
    \caption{Illustration of a half-edge $h$, $\tgt(h)$, $\opp(h)$, and $\lft(h)$.}
    \label{fig:basicDefn}
\end{figure}
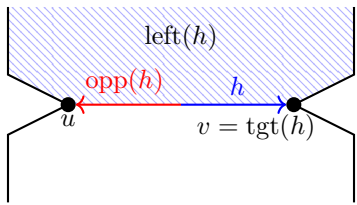

We view each edge of the graph as consisting of
two oppositely directed half-edges, see Figure~\ref{fig:basicDefn} for illustration.  For a half-edge $h$, let $\tgt(h)$ denote the vertex of $G$ towards which $h$ is directed,
let $\opp(h)$ denote the opposite half-edge, and let $\lft(h)$ denote the face of $G$ drawn to the left of $h$.  Let $\hes(G)$
denote the set of half-edges of $G$.  It will be often convenient to select just one half-edge from each pair arbitrarily; let $\ohes(G)$
denote a subset of $\hes(G)$ containing exactly one half-edge from each pair of oppositely directed half-edges.

It would be needlessly complicated for algorithms to operate on drawings as defined at the beginning of the section (we would need
to come up with a discrete way to describe the curves).  For the purposes of our algorithm, a graph $G$ with a 2-cell drawing is represented
by giving
\begin{itemize}
\item the set $V(G)$ of vertices of $G$, the set $\hes(G)$ of half-edges of $G$, and the set $F(G)$ of faces of $G$, where
the faces are taken as abstract elements rather than subsets of the surface, and
\item the functions $\opp:\hes(G)\to \hes(G)$, $\tgt:\hes(G)\to V(G)$, and $\lft:\hes(G)\to F(G)$.
\end{itemize}
We also occassionally refer to the set $E(G)$ of edges of $G$, which can be viewed as pairs $\{h,\opp(h)\}$ for $h\in \hes(G)$.
Let us remark that information contained in this representation determines the drawing of $G$ up to homeomorphisms of $\Sigma$.
The \emph{size} of $G$ is defined as $|V(G)|+|F(G)|+|E(G)|$.  Note that a graph and its dual have the same size,
and generalized Euler's formula implies that if $G$ is drawn in a surface
of Euler genus $g$ and $G^\star$ is a simple graph with $n$ vertices, then the size of $G$ is $O(n+g)$.

\subsection*{Homology}

We are going to need some simple definitions from the homology theory.
\begin{itemize}
\item A \emph{$0$-chain} is a formal sum of vertices of $G$ with integer coefficients; $0$-chains form a free abelian group $C_0(G)$.
For a $0$-chain $b$ and a vertex $v$, let $b[v]$ denote the coefficient at $v$, and let $|b|=\sum_{v\in V(G)} |b[v]|$.
\item A \emph{$1$-chain} is a formal sum $K=\sum_{h\in \ohes(G)} c_h\cdot h$ with integer coefficients,
i.e., $1$-chains form a free abelian group $C_1(G)$ generated by $\ohes(G)$.
For a half-edge $h$, let us define $K[h]=c_h$ if $h\in\ohes(G)$ and $K[h]=-c_{\opp(h)}$ otherwise.
Note that 
\begin{align}
K[\opp(h)]=-K[h].
\label{eq:opp}
\end{align}  
We view each half-edge $h\in \hes(G)\setminus\ohes(G)$ as the $1$-chain $-\opp(h)$.
For example, a directed walk $W$ traversing half-edges $h_1$, \ldots, $h_m$ corresponds to the $1$-chain $K=h_1+\cdots+h_m$,
and for each half-edge $h$, $K[h]$ is the number of times $W$ traverses $h$ minus the number of times it traverses $\opp(h)$.
Let us also define $|K|=\sum_{h\in\ohes(G)} |K[h]|$.
\item A \emph{$2$-chain} is a formal sum of faces of $G$ with integer coefficients.  Let $C_2(G)$ be the free abelian group of $2$-chains.
For a $2$-chain $A$ and a face $x$, let $A[x]$ denote the coefficient of $x$ in $A$.
\end{itemize}

\begin{figure}
\centering{
\includegraphics[page=3]{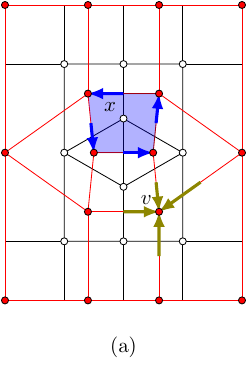}
\hskip 2em
\includegraphics[page=7]{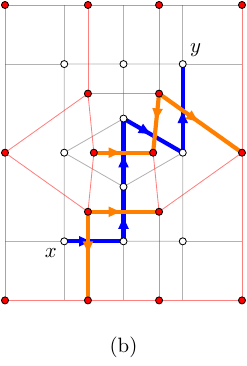}
}
\caption{(a) Example of $\partial_2 x$ and $\partial_2^\star v$ in a graph $G$ (red) for a face $x$ and a vertex $v$ of $G$.
(b) A path from face $x$ to face $y$ in the dual (blue) and the corresponding copath in $G$ (orange).}
\label{fig:boundary}
\end{figure}

Let us now define boundary operators:
\begin{itemize}
\item For a face $x\in F(G)$, let $$\partial_2 x=\sum_{h\in \hes(G):\lft(h)=x} h.$$
That is, the boundary of $x$ consists of the incident half-edges that are directed counter-clockwise around $x$, see Figure~\ref{fig:boundary}(a) for an example.
Let us extend $\partial_2$ to $2$-chains linearly.
\item To each half-edge $h$, we assign a $0$-chain $\partial_1 h=\tgt(h)-\tgt(\opp(h))$,
and we extend the mapping $\partial_1$ to all $1$-chains linearly.  For example, if $W$ is the 1-chain corresponding to
a walk starting in a vertex $u$ and ending in a vertex $v$, then $\partial_1 W=v-u$.
\item For each $0$-chain $b=\sum_{v\in V(G)} b_v\cdot v$, we define $\partial_0 b=\sum_{v\in V(G)} b_v$.
\end{itemize}

\noindent Now, for $i\in\{0,1\}$,
\begin{itemize}
\item an \emph{$i$-boundary} is an $i$-chain belonging to the subgroup $B_i(G)=\{\partial_{i+1} a:a\in C_{i+1}(G)\}$ of $C_i(G)$, and
\item an \emph{$i$-cycle} is an $i$-chain belonging to the subgroup $Z_i(G)=\{K\in C_i(G):\partial_i K=0\}$.
For example, if $W$ is the 1-chain corresponding to a closed walk, then $W$ is a $1$-cycle.
\end{itemize}
Note that
\begin{itemize}
\item every $i$-boundary is an $i$-cycle, and thus $B_i(G)$ is a subgroup of $Z_i(G)$;
\item since $G$ is connected, every $0$-cycle is a $0$-boundary, and thus $B_0(G) = Z_0(G)$; and
\item for any contractible cycle $C$ in $G$, the corresponding $1$-chain is a $1$-boundary, since it can be expressed
as the sum of $\partial_2 x$ over the faces $x$ drawn in the open disk bounded by $C$.
\end{itemize}

The \emph{(first) homology group} $H_1(G)$ is defined as the quotient $Z_1(G)/B_1(G)$.
For example, suppose that $W_1$ and $W_2$ are closed walks and $W_2$ is freely homotopic to $W_1$.  Then $W_2$ can be obtained from $W_1$
by a sequence of the operations of adding a walk around a face and removing the subwalks consisting of taking an edge in one direction
and immediately coming back over it.  Hence, for the corresponding $1$-chains, $W_2=W_1+Q$ for some $1$-boundary $Q$, and thus
$W_1$ and $W_2$ correspond to the same element of $H_1(G)$.

Let us also define the dual (cohomology) operators:
\begin{itemize}
\item For a vertex $v$, let $\partial^\star_2 v=\sum_{h\in \hes(G):\tgt(h)=v} h$ and let us extend $\partial^\star_2$ to $0$-chains
linearly, see Figure~\ref{fig:boundary}(a) for an example.
\item For a half-edge $h$, let $\partial^\star_1 h=\lft(h)-\lft(\opp(h))$ and let us extend $\partial^\star_1$ to all $1$-chains
linearly.
\item For a $2$-chain $d=\sum_{x\in F(G)} d_x\cdot x$, let $\partial^\star_0 d=\sum_{x\in F(G)} d_x$.
\end{itemize}
\noindent And,
\begin{itemize}
\item a \emph{coboundary} is a $1$-chain belonging to the subgroup $B^\star(G)=\{\partial^\star_2 a:a\in C_0(G)\}$ of $C_1(G)$,
\item a \emph{cocycle} is a $1$-chain belonging to the subgroup $Z^\star(G)=\{K\in C_1(G):\partial^\star_1 K=0\}$, and
\item the \emph{cohomology group} $H^\star(G)=Z^\star(G)/B^\star(G)$.
\end{itemize}
Note that coboundaries correspond to edge cuts in $G$, similarly to the way 1-boundaries correspond to separating cycles.
The \emph{dual graph} $G^\star$ to $G$ is a graph with vertex set $F(G)$ and with each edge $e$ of $G$ corresponding to an
edge $e^\star$ of $G^\star$ joining the faces incident with $e$.  Note that rather than defining the dual boundary operators,
we could work with chains in the dual graph; however, having both primal and dual boundary operators act on the same sets
simplifies the notation.  For two faces $x$ and $y$ of $G$, a \emph{copath from $x$ to $y$} is a $1$-chain $P$ with $\partial^\star_1 P=y-x$,
see Figure~\ref{fig:boundary}(b) for an illustration.
Let us remark that a cocycle is a copath from any face to itself.

For a $1$-chain $f$ in $G$ and $h\in \hes(G)$, define $\trans{f}(h)=f[h]$,
and let us extend the function $\trans{f}$ to all $1$-chains $K$ in $G$ linearly. The following observation based on \eqref{eq:opp} will be helpful in future calculations.
\begin{observation}
Let $G$ be a graph with a $2$-cell drawing in an orientable surface and let $f$ and $K$ be $1$-chains in $G$. Then
\[\trans{f}(K) =
\sum_{h\in \ohes(G)} f[h]K[h]
=
\sum_{h \in \hes(G): K[h]>0} f[h]K[h].\]
\end{observation}
In a typical application, we view $f$ as a flow with excess $(\partial_1 f)[v]$ at each vertex $v$
(see the next section for details) and $K$ as a cocycle.  In this case, $\trans{f}(K)$ gives the amount of flow $f$ sends over $K$.
In particular, as one would expect, the following relation holds in the case that $K$ consists of half-edges entering a vertex.
\begin{observation}\label{obs-transv}
Let $G$ be a graph with a 2-cell drawing in an orientable surface, let $f$ be a $1$-chain in $G$, let $v$ be a vertex of $G$,
and let $K=\partial^\star_2 v$.  Then
$$\trans{f}(K)=(\partial_1 f)[v].$$
Moreover, if $f$ is a $1$-cycle, then $\trans{f}(R)=0$ for every coboundary $R$.
\end{observation}
\begin{proof}
Since $K=\partial^\star_2 v=\sum_{h\in \hes(G):\tgt(h)=v} h$, we have
\begin{align*}
\trans{f}(K)&=\sum_{h\in \hes(G):\tgt(h)=v} \trans{f}(h)=\sum_{h\in\hes(G):\tgt(h)=v} f[h]\\
&=\sum_{h\in\ohes(G)} f[h]\cdot (\partial_1 h)[v]=(\partial_1 f)[v].
\end{align*}
If $f$ is a $1$-cycle, i.e., $\partial_1 f=0$, then this implies $\trans{f}(\partial^\star_2 v)=0$ for each vertex $v$, and
we can extend this claim linearly to all coboundaries.
\end{proof}

Let us note the dual form of this observation.
\begin{observation}\label{obs-transf}
Let $G$ be a graph with a 2-cell drawing in an orientable surface, let $x$ be a face of $G$, let $K$ be a $1$-chain in $G$,
and let $f=\partial_2 x$.  Then
$$\trans{f}(K)=(\partial^\star_1 K)[x].$$
Moreover, if $K$ is a cocycle, then $\trans{b}(K)=0$ for every 1-boundary $b$.
\end{observation}

Let us also note the following well-known fact, describing how to obtain a basis of the first homology group of a graph drawn
in an orientable surface.

\begin{figure}
\centering{
\includegraphics[page=4,scale=0.7]{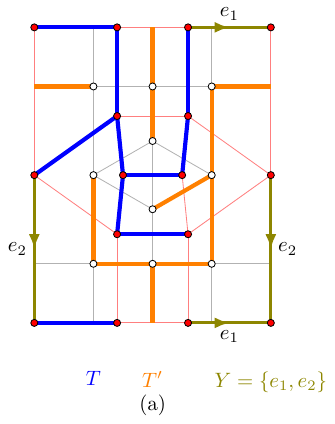}
\includegraphics[page=5,scale=0.7]{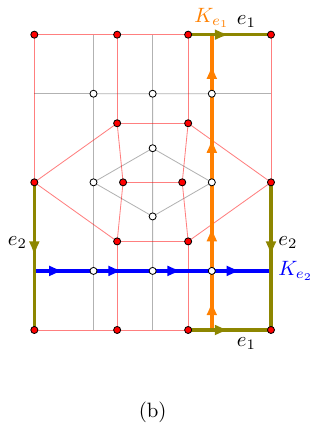}
\includegraphics[page=6,scale=0.7]{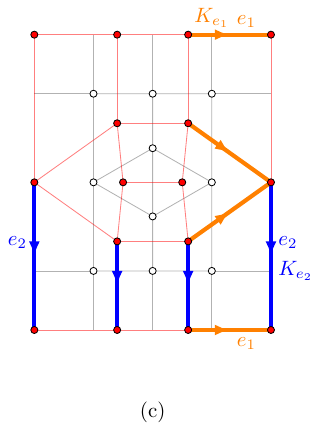}
\includegraphics[page=13,scale=0.7]{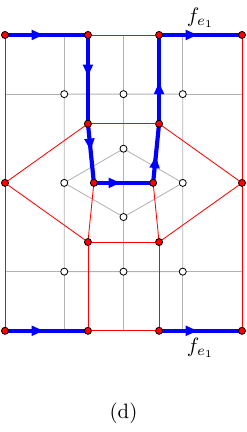}
\includegraphics[page=14,scale=0.7]{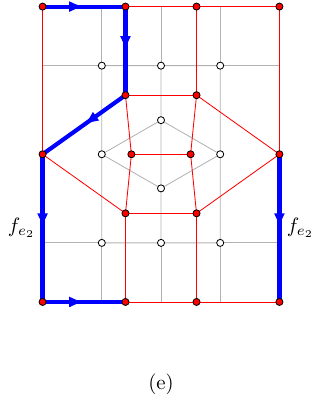}
}
\caption{Example of obtaining a basis of the first homology group as described in Observation~\ref{obs-genhom}.
(a) Depicts a spanning tree $T$ of $G$ (blue), a spanning tree $T'$ of $G^\star-X$ (orange), and the set $Y=\{e_1,e_2\}$ (green),
(b) depicts cycles in $T'+e_1^\star$ (orange) and $T'+e_2^\star$ (blue),
(c) depicts the corresponding cocycles $K_{e_1}$ and $K_{e_2}$ forming
the basis of $H^\star(G)$, (d) and (e) depict the cycles $f_{e_1}$ and $f_{e_2}$
in $T+e_1$ and $T+e_2$ forming the basis of $H_1(G)$. }
\label{fig:cocycleBasis}
\end{figure}

\begin{observation}\label{obs-genhom}
Let $G$ be a graph with a 2-cell drawing in an orientable surface $\Sigma$ of Euler genus $g$.
Both the homology group and the cohomology group of $G$ are isomorphic to $\mathbb{Z}^{g}$.
Moreover, there is an algorithm that in time linear in the size of $G$ returns their bases
$M=\{f_e:e\in Y\}$ and $Q=\{K_e:e\in Y\}$ indexed by a set $Y$ of size $g$, such that
for $e,e'\in Y$, we have
\begin{equation}
\label{eq:feKe}
\trans{f_e}(K_{e'})=\begin{cases}
1&\text{ if $e=e'$}\\
0&\text{ otherwise.}
\end{cases}
\end{equation}
\end{observation}
\begin{proofout}
The bases can be obtained as follows (see Figure~\ref{fig:cocycleBasis} for an illustration): Let $T$ be a spanning tree of $G$ and let $X=\{e^\star:e\in E(T)\}$.
Observe that $G^\star-X$ is connected, let $T'$ be a spanning tree of $G^\star-X$,
and let $S$ be the subgraph of $G$ with $V(S)=V(G)$ and $E(S)=\{e:e^\star\in E(T')\}$.
Observe that the set $Y=E(G)\setminus (E(T)\cup E(S))$ has size $g$.
For each $e\in Y$, choose one of the half-edges $h_e$ forming $e$, let $f_e$ be the unique $1$-cycle in $C_1(T+e)$ such that $f_e[h_e]=1$
and let $K_e$ be the unique cocycle in $C_1(S+e)$ such that $K_e[h_e]=1$.
Then $M=\{f_e:e\in Y\}$ is a basis of $H_1(G)$ and $Q=\{K_e:e\in Y\}$ is a basis of $H^\star(G)$ satisfying (\ref{eq:feKe}).
\end{proofout}

\begin{corollary}\label{cor-testb}
Let $G$ be a graph with a 2-cell drawing in an orientable surface.
A $1$-cycle $f$ is a $1$-boundary if and only if $\trans{f}(K)=0$
for every $K\in H^\star(G)$.
\end{corollary}
\begin{proof}
By Observation~\ref{obs-transf}, $\trans{f}(K)=0$ for every $1$-boundary $f$ and cocycle $K$.  Conversely,
suppose that $\trans{f}(K)=0$ for every $K\in H^\star(G)$.
Let $Y$, $M=\{f_e:e\in Y\}$ and $Q=\{K_e:e\in Y\}$ be as in Observation~\ref{obs-genhom}.
Since $f$ is a $1$-cycle and $M$ is a basis of $H_1(G)$, we have $f=b+\sum_{e\in Y}\alpha_e f_e$ for some $1$-boundary
$b$ and integers $\alpha_e$ for $e\in Y$.  For each $e'\in Y$, linearity and (\ref{eq:feKe}) give
$$0=\trans{f}(K_{e'})=\trans{b}(K_{e'})+\sum_{e\in Y} \alpha_e\trans{f_e}(K_{e'})=\alpha_{e'},$$
and thus $f=b$ is a $1$-boundary.
\end{proof}

\subsection*{Flow-coloring duality}

We say that a $1$-chain $f$ is a \emph{flow} if $f[h]\in\{-1,0,1\}$
for each $h\in \hes(G)$---we view $f$ as sending the amount $f[h]$ in the direction of $h$.
Note that for each $v\in V(G)$, $(\partial_1 f)[v]$ is the excess of the flow $f$ in $v$;
we consider flows with sources and sinks, and thus, in general, we do not require that $\partial_1 f=0$.
We say that a flow $f$ is \emph{nowhere-zero} if $f[h]\neq 0$ for each $h\in \hes(G)$.
We say that a $0$-chain $d$ is \emph{divisible by $k$} if $k|d[v]$ for each $v\in V(G)$.

\begin{figure}
\begin{center}
\includegraphics[page=1]{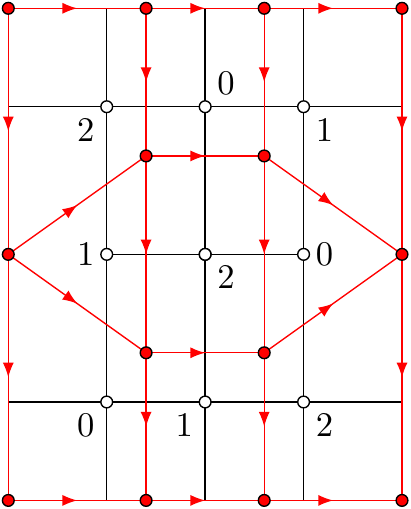}
\end{center}
\caption{A homomorphism of $G^\star$ to a triangle and the corresponding nowhere-zero flow $f$ in $G$.}\label{fig-homomorphism}
\end{figure}

The duality between flows and colorings was discovered by Tutte~\cite{tutteflow},
and the version for homomorphisms to cycles by Goddyn et al.~\cite{goddyn}.
We include the proof to account for the differences in terminology.
We define the vertex set of a cycle $C$ of length $m$ to be $\{0,\ldots,m-1\}$, with
each vertex $v$ adjacent to $v-1$ and $v+1\pmod m$.
Figure~\ref{fig-homomorphism} gives an example of a homomorphism to a triangle
and the corresponding nowhere-zero flow in the dual graph.

\begin{lemma}\label{lemma-tutte}
Let $G$ be a graph with a 2-cell drawing in an orientable surface, let $S$ be a non-empty subset of $F(G)$,
let $C$ be a cycle of length $m$, and let $\psi:S\to V(C)$ be an arbitrary function.  Let $x$ be an arbitrary element of $S$ and
for each $y\in S\setminus\{x\}$, let $P_y$ be a copath from $x$ to $y$.
Let $Q$ be a basis of $H^\star(G)$.  Then the following claims are equivalent:
\begin{itemize}
\item[(i)] $\psi$ extends to a homomorphism from $G^\star$ to $C$.
\item[(ii)] There exists a nowhere-zero flow $f$ such that
\begin{itemize}
\item $\partial_1 f$ is divisible by $m$,
\item $\trans{f}(K)$ is divisible by $m$ for every $K\in Q$, and
\item for each $y\in S\setminus\{x\}$, we have $\trans{f}(P_y)\equiv \psi(y)-\psi(x)\pmod m$.
\end{itemize}
\end{itemize}
\end{lemma}
\begin{proof}
Suppose $\varphi:F(G)\to V(C)$ is a homomorphism from $G^\star$ extending $\psi$,
and let us extend $\varphi$ to $2$-chains in $G$ linearly.
Let $f$ be the nowhere-zero flow in $G$ such that for each $h\in \hes(G)$, $f[h]$ is the unique element of $\{-1,1\}$
satisfying
$$f[h]\equiv \varphi(\partial^\star_1 h)=\varphi(\lft(h))-\varphi(\lft(\opp(h)))\pmod m;$$
such an element exists since $\varphi$ is a homomorphism to $C$.
For any $1$-chain $K$ we have
$$\trans{f}(K)=\sum_{h\in\ohes(G)} K[h]f[h]\equiv \sum_{h\in\ohes(G)} K[h]\varphi(\partial^\star_1 h)=\varphi(\partial^\star_1 K)\pmod m$$
by linearity. In particular, if $K$ is a cocycle, i.e., $\partial^\star_1 K=0$, then $\trans{f}(K)\equiv 0\pmod m$. Therefore,
\begin{itemize}
\item $\trans{f}(K)\equiv 0\pmod m$ for every $K\in Q$, and
\item by Observation~\ref{obs-transv}, for every vertex $v\in V(G)$ we have
$$(\partial_1 f)[v]=\trans{f}(\partial^\star_2 v)\equiv 0\pmod m,$$
and thus $\partial_1 f$ is divisible by $m$.
\end{itemize}
Finally, for each $y\in S\setminus\{x\}$, we have
\begin{align*}
\trans{f}(P_y)&\equiv \varphi(\partial^\star_1 P_y)=\varphi(y-x)\\
&=\varphi(y)-\varphi(x)=\psi(y)-\psi(x)\pmod m.
\end{align*}
Hence, (i) implies (ii).

Conversely, suppose that $f$ is a nowhere-zero flow satisfying the conditions listed in (ii).
Let $Q'=\{\partial^\star_2 v:v\in V(G)\}$.
By the assumptions and Observation~\ref{obs-transv}, we have $m|\trans{f}(K_0)$ for each $K_0\in Q\cup Q'$.
Since $Q$ generates $H^\star(G)$ and $Q'$ generates $B^\star(G)$, $Q\cup Q'$ generates
$Z^\star(G)$, and thus  $m|\trans{f}(K)$ for each cocycle $K$.
For each $y\in F(G)\setminus (S\setminus\{x\})$, let $P_y$ be an arbitrary copath from $x$ to $y$.
For each $y\in F(G)$, we define
$$\varphi(y)=(\psi(x)+\trans{f}(P_y)) \bmod m.$$ 
We claim that $\varphi$ is a homomorphism from $G^\star$ to $C$ extending $\psi$.
Indeed, suppose that $y,y'\in F(G)$ are adjacent vertices of $G^\star$, and
let $h$ be a half-edge of $G$ with $\lft(h)=y$ and $\lft(\opp(h))=y'$.  Then
$$\partial^\star_1 (P_y-h-P_{y'})=(y-x)-(y-y')-(y'-x)=0,$$
and thus $P_y-h-P_{y'}$ is a cocycle and $m|\trans{f}(P_y-h-P_{y'})$.  Hence,
\begin{align*}
\varphi(y)-\varphi(y')&\equiv \trans{f}(P_y)-\trans{f}(P_{y'})\\
&=\trans{f}(h)+\trans{f}(P_y-h-P_{y'})\\
&\equiv f[h]=\pm1\pmod m,
\end{align*}
since $f$ is nowhere-zero.  It follows that $\varphi$ is
a homomorphism from $G^\star$ to $C$.
Moreover, for any $y\in S$, we have
$$\varphi(y)\equiv \psi(x)+\trans{f}(P_y)\equiv \psi(y)\pmod m,$$
if $y\neq x$ by the assumptions and if $y =x$ since $P_x$ is a cocycle.
Therefore $\varphi(y)=\psi(y)$ for each $y\in S$, and thus (ii) implies (i).
\end{proof}

For $1$-chains $K$ and $K'$, we write $K\preceq K'$ if $0\le K[h]\le K'[h]$ for
every half-edge $h$ such that $K'[h]\ge 0$.
For a 1-chain $f$ in $G$, an \emph{$f$-circulation} is a 1-cycle $c\preceq f$ in $G$; in particular,
$c$ can only send flow in the same direction as $f$ does.
Circulations can be used to translate between any two nowhere-zero flows with the same boundary.
\begin{observation}\label{obs-reverse}
Suppose $f_0$ is a nowhere-zero flow in a graph $G$ and let $d=\partial_1 f_0$.
\begin{itemize}
\item If $c$ is an $f_0$-circulation, then $f_0-2c$ is a nowhere-zero flow in $G$
satisfying $\partial_1 (f_0-2c)=d$.
\item If $f$ is a nowhere-zero flow in $G$ satisfying $\partial_1 f=d$, then $(f_0-f)/2$ is an $f_0$-circulation.
\end{itemize}
\end{observation}

\begin{corollary}\label{cor-circ}
Let $G$ be a graph, let $d$ be a $0$-boundary, and let $f_0$ be a nowhere-zero flow in $G$ such that $\partial_1 f_0=d$.
Let $Q$ be a system of pairs $(K,a)$, where $K$ is a $1$-chain and $a$ is an integer.
The following claims are equivalent:
\begin{itemize}
\item There exists a nowhere-zero flow $f$ in $G$ such that $\partial_1 f=d$ and $\trans{f}(K)=a$ for each $(K,a)\in Q$.
\item There exists an $f_0$-circulation $c$ such that $\trans{c}(K)=\frac{1}{2}(\trans{f_0}(K)-a)$ for each $(K,a)\in Q$.
\end{itemize}
\end{corollary}

Note that for a fixed $0$-boundary $d$, Corollary~\ref{cor-circ} splits the verification of the condition (ii) of Lemma~\ref{lemma-tutte} to two independent
steps.
\begin{itemize}
\item[(A)] Deciding whether there exists a nowhere-zero flow $f_0$ in $G$ with $\partial_1 f_0=d$.
\item[(B)] Deciding whether there exists an $f_0$-circulation $c$ with the prescribed values of $\trans{c}$
over fixed copaths and cocycles (modulo $m$).
\end{itemize}

The first part (A) is easily achieved using any efficient maximum flow algorithm in polynomial time,
as shown in the following lemmas.  We say that a $0$-boundary $d$ in a graph $G$ is \emph{parity-compliant}
if for each $v\in V(G)$, $d[v]$ and the degree of $v$ in $G$ have the same parity.

\begin{lemma}\label{lemma-flow1}
Let $G$ be a graph and let $d$ be a $0$-boundary.  The following claims are equivalent.
\begin{itemize}
\item[(i)] There exists a nowhere-zero flow $f_0$ in $G$ with $\partial_1 f_0=d$.
\item[(ii)] The $0$-boundary $d$ is parity-compliant and there exists a flow $f_1$ in $G$ with $\partial_1 f_1=d$.
\end{itemize}
Moreover, given $f_1$ as in (ii), a nowhere-zero flow $f_0$ with $\partial_1 f_0=d$ can be found in linear time.
\end{lemma}
\begin{proof}
Suppose $f_0$ exists; since $|f_0[h]|=1$ for each half-edge $h$, $d[v]=(\partial_1 f_0)[v]=\trans{f_0}(\partial^\star_2 v)$ has the same parity as the degree of $v$ for each $v\in V(G)$.
Hence, $d$ is parity-compliant, and we can set $f_1=f_0$.

Conversely, suppose that (ii) holds.  Consider the undirected graph $T$ with vertex set $V(G)$ and the edge set consisting of the edges $e\in E(G)$
such that $f_1[h]=0$ for either (or equivalently, both) of the half-edges $h$ of $e$.  Since $d=\partial_1 f_1$ is
parity-compliant, every vertex has even degree in $T$, and thus in every component of $T$, there exists a closed walk passing through
every edge of the component exactly once.  For each such walk, choose a direction, and let $f'_1\in C_1(G)$ be the sum of the half-edges
of $T$ whose direction matches the direction selected for the walk that contains it.  Then $f'_1$ is a flow in $G$ and $\partial_1 f'_1=0$,
and for each half-edge $h$, exactly one of $|f_1[h]|$ and $|f'_1[h]|$ is one and the other one is zero.
Consequently, $f_0=f_1+f'_1$ is a nowhere-zero flow in $G$ and $\partial_1 f_0=\partial_1 f_1=d$.  Moreover, the flow $f'_1$
can be found in linear time using the standard algorithm to find Eulerian tours.
\end{proof}

Hence, (A) reduces to finding any flow with the given parity-compliant boundary $d$. As is well known, this can
be restated in terms of finding a flow in a corresponding \emph{network} (directed graph with edges of bounded capacities
and with two vertices designated as a source and a sink).

\begin{lemma}\label{lemma-flow2}
Let $G$ be a graph and let $d$ be a $0$-boundary.  Let $H$ be the network obtained as follows: Start with $G$
and replace each edge by a pair of oppositely directed edges with capacity $1$.  Add a vertex $s$ and for each $v\in V(G)$
such that $d[v]<0$, add an edge from $s$ to $v$ of capacity $|d[v]|$; and add a vertex $t$ and for each $v\in V(G)$ such that $d[v]>0$, add an edge from $v$ to $t$ of capacity $d[v]$.
There exists a flow $f_1$ in $G$ with $\partial_1 f_1=d$ if and only if the network $H$ contains a flow of size $|d|/2$ from $s$ to $t$.
\end{lemma}
\begin{proof}
Note that $f_1$ can be turned into a flow in $H$ of size $|d|/2$ by setting the flow over each edge between $v\in V(G)$ and $s$ or $t$ to $|d[v]|$.

Conversely, since $H$ has integer capacities, we can assume that $H$ contains a network flow $g$ of size $|d|/2$ from $s$ to $t$
with integer values.  Since $g$ has size $|d|/2$, the flow over each edge between $v\in V(G)$ and $s$ or $t$ is equal to $|d[v]|$.  Hence,
interpreting the half-edges of $G$ as the directed edges of $H-\{s,t\}$ in the natural way and defining
$f_1=\sum_{h\in \ohes(G)} (g(h)-g(\opp(h))\cdot h$,
we conclude that $f$ is a flow satisfying $\partial_1 f_1=d$.
\end{proof}

Finally, let us note a bound on the number of possible boundaries that we are going to need to test.
We say that a $0$-boundary $d$ in a graph $G$ is \emph{relevant} if $d$ is parity-compliant and $|d[v]|$ is at most
the degree of $v$ for each $v\in V(G)$.
\begin{observation}\label{obs-relev}
Let $G$ be a graph drawn in an orientable surface and let $H=G^\star$ be its dual.
For any nowhere-zero flow $f$ in $G$, the $0$-boundary $\partial_1 f$ is relevant.
Moreover, if $C$ is a cycle, then there are at most $q^\star_C(H)$ relevant $0$-boundaries in $G$ divisible by $|C|$,
and each such $0$-boundary $d$ satisfies $|d|+1\le b^\star_C(H)$.
\end{observation}

To solve the second part (B), we use the algorithms that we describe in the next section.  

\section{Circulations with prescribed homology}\label{sec-homol}

For $1$-chains $f$ and $K$, we define
$$\transplus{f}(K)=\sum_{h\in \hes(G): f[h],K[h]>0} f[h]K[h].$$
Intuitively, if $K$ is a cocycle, then $\transplus{f}(K)$ is an
upper bound on the amount of flow that can be sent over $K$ by
any flow $f'$ such that $f'\preceq f$.
Let us remark on the following basic properties of this notion.
\begin{observation}\label{obs-rev}
For any $1$-chains $f$ and $K$, we have
$$\transplus{f}(K)=\transplus{-f}(-K)$$
and
$$\trans{f}(K)=\transplus{f}(K)-\transplus{-f}(K),$$
and if $f$ is a nowhere-zero flow, then
$$|K|=\transplus{f}(K)+\transplus{-f}(K).$$
\end{observation}
\begin{proof}
The first claim holds since for every $h\in \hes(G)$, we have
$(-f)[\opp(h)]=f[h]$, $(-K)[\opp(h)]=K[h]$, and $\opp$ is an involution.
The second claim holds since
\begin{align*}
    \trans{f}(K)&=\sum_{h\in \hes(G):K[h]>0} f[h]K[h]\\
    &=\sum_{h\in \hes(G):K[h]>0} (\max(f[h],0)-\max(-f[h],0))K[h]=\transplus{f}(K)-\transplus{-f}(K).
\end{align*}
For the last claim, we assume that $f$ is a nowhere-zero flow, i.e., $|f[h]|=1$ for every half-edge $h$.
Consequently,
\begin{align*}
|K|&=\sum_{h\in \hes(G):K[h]>0} K[h]\\
&=\sum_{h\in \hes(G):K[h]>0} (\max(f[h],0)+\max(-f[h],0))K[h]=\transplus{f}(K)+\transplus{-f}(K).
\end{align*}
\end{proof}

For a $1$-chain $K$ in a graph $G$, let $\transplus{f}(K+B^\star(G))$ denote the minimum of $\transplus{f}(K+R)$ over all coboundaries~$R$.
The motivation for this definition is the following necessary condition on $f$-circulations.
\begin{lemma}\label{lemma-nomore}
Let $f$ and $K$ be $1$-chains in a graph $G$.
For every $f$-circulation $c$,
$$\trans{c}(K)\le \transplus{f}(K+B^\star(G)).$$
\end{lemma}
\begin{proof}
Let $R$ be a coboundary such that letting $K'=K+R$, we have
$\transplus{f}(K+B^\star(G))=\transplus{f}(K')$.
By Observation~\ref{obs-transv} we have $\trans{c}(R)=0$.
Since $c$ is an $f$-circulation, we have
\begin{align*}
\trans{c}(K)&=\trans{c}(K')-\trans{c}(R)=\trans{c}(K')\\
&=\sum_{h\in \hes(G):c[h]>0} c[h]K'[h]\le \sum_{h\in \hes(G):c[h],K'[h]>0} c[h]K'[h]\\
&\le \sum_{h\in \hes(G):f[h],K'[h]>0} f[h]K'[h]=\transplus{f}(K')=\transplus{f}(K+B^\star(G)).
\end{align*}
\end{proof}

In (B), we ask for an $f_0$-circulation $c$ with the prescribed values of $\trans{c}$
over fixed copaths and cocycles.  We start by showing that there at least always exists a 1-cycle $b$
with the prescribed values of $\trans{b}$ over the relevant copaths and cocycles.
This is a straightforward consequence of the orthogonality of the bases of $H_1(G)$ and $H^\star(G)$.

\begin{figure}
\begin{center}
\includegraphics[page=8]{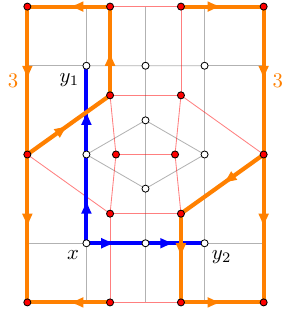}
\end{center}
\caption{
Example for Lemma~\ref{lemma-extflow} with
$S=\{x,y_1,y_2\}$, $a'(x)=0$, $a'(y_1)=1$, $a'(y_2)=1$,
$a(K_{e_1})=0$, $a(K_{e_2}) = 2$; see Figure~\ref{fig:cocycleBasis} for the choice of the bases.
Copaths $P_{y_1}$ and $P_{y_2}$ are depicted in $G^\star$ as blue. The desired 1-cycle $b$ is orange, and it is
obtained as $b = 2f_{e_2} + 3 \partial_2 y_1 + \partial_2 y_2$. One edge of $b$ has coefficient 3, the others are with coefficient 1.}\label{fig:extflow}
\end{figure}

\begin{lemma}\label{lemma-extflow}
Let $G$ be a graph of size $n$ with a 2-cell drawing in an orientable surface of Euler genus $g$.
Let $a:H^\star(G)\to\mathbb{Z}$ be a homomorphism.
Let $S$ be a non-empty subset of $F(G)$, let $x$ be an element of $S$,
for each $y\in S$ let $P_y$ be a copath from $x$ to $y$,
and let $a':S\to \mathbb{Z}$ be an arbitrary function such that $a'(x)=0$.
Then there exists a $1$-cycle $b$ such that $\trans{b}(K)=a(K)$ for each $K\in H^\star(G)$
and $\trans{b}(P_y)=a'(y)$ for each $y\in S$.  Moreover, $b$ can be constructed
in time $O((g+1)|S|n)$.
\end{lemma}
\begin{proof}
Let $Y$, $M=\{f_e:e\in Y\}$ and $Q=\{K_e:e\in Y\}$ be as in Observation~\ref{obs-genhom}.
The 1-cycle $b$ is obtained from a linear combination of the elements of the basis $M$
chosen so that $\trans{b}(K_e)=a(K_e)$ for each $e\in Y$, and thus also $\trans{b}(K)=a(K)$
for each $K\in H^\star(G)$. To ensure that $\trans{b}(P_y)=a'(y)$ for each
$y\in S$, we add a suitable multiple of $\partial_2 y$ to $b$;
by Observation~\ref{obs-transf}, this does not affect the values of $\trans{b}(K)$ for $K\in H^\star(G)$.
See Figure~\ref{fig:extflow} for an example.

More precisely, for each $y\in S\setminus\{x\}$, let $\gamma_y \in  \mathbb{Z}$ be
\begin{align}\label{eq:gamma}
    \gamma_y=a'(y)-\sum_{e\in Y} a(K_e)\trans{f_e}(P_y).
\end{align}
Let 1-cycle $b$ be
\begin{align}\label{eq:b}
    b=\sum_{e\in Y} a(K_e)f_e + \sum_{y\in S\setminus\{x\}} \gamma_y\partial_2 y.
\end{align}
Consider any cocycle $K\in H^\star(G)$, where $K=\sum_{e\in Y}\alpha_eK_e$ for some integers $\alpha_e$.
Note that for $y\in F(G)$, we have $\trans{\partial_2 y}(K)=0$ by Observation~\ref{obs-transf}, and by \eqref{eq:feKe},
it follows that
\begin{align*}
\trans{b}(K)&=\sum_{e\in Y} a(K_e)\sum_{e'\in Y}\alpha_{e'}\trans{f_e}(K_{e'})=\sum_{e\in Y} a(K_e)\alpha_e\\
&=a\left(\sum_{e\in Y} \alpha_e K_e\right)=a(K),
\end{align*}
since $a$ is a homomorphism.  Moreover, for each $y\in S\setminus\{x\}$ and a face $y'\not\in \{x,y\}$, we have
$\trans{\partial_2 y'}(P_y)=(\partial^\star_1 P_y)[y']=0$, and 
$\trans{\partial_2 y}(P_y)=(\partial^\star_1 P_y)[y]=1$, implying that
$$\trans{b}(P_y)=\left(\sum_{e\in Y} a(K_e)\trans{f_e}(P_y)\right)+\gamma_y=a'(y)$$
by the choice of $\gamma_y$.

Calculating $\gamma_y$ in \eqref{eq:gamma} takes time $O(1+gn)$. Calculating
$b$ from \eqref{eq:b} then takes time $O(gn + |S|((gn+1)+n)) = O((g+1)|S|n)$.
\end{proof}

Let us remark that for the purposes of the algorithm mentioned in the statement of
Lemma~\ref{lemma-extflow}, the homomorphism $a$ should be given by its values on
the basis $Q$ of $H^\star(G)$ obtained using Observation~\ref{obs-genhom}.
This is of course not a substantial restriction, as if it were given in any other basis,
we could just transform it; and additionally, in all the uses $a$ is in fact represented in this way.

The proof of the following theorem, which is the cornerstone of our approach, is inspired by Chambers et al.~\cite{homologycuts,homologycutsB}.
A $1$-chain $K$ is \emph{simple} if $|K[h]|\le 1$ for every half-edge $h$.

\begin{theorem}\label{thm-extcirc}
Let $G$ be a graph of size $n$ with a 2-cell drawing in an orientable surface of Euler genus $g$.
Let $a:H^\star(G)\to\mathbb{Z}$ be a homomorphism.
Let $S$ be a non-empty subset of $F(G)$, let $x$ be an element of $S$,
for each $y\in S$ let $P_y$ be a copath from $x$ to $y$,
and let $a':S\to \mathbb{Z}$ be an arbitrary function such that $a'(x)=0$.
Let $f$ be a 1-chain in $G$.  The following claims are equivalent:
\begin{itemize}
\item[(i)] There exists an $f$-circulation $c$ such that $\trans{c}(K)=a(K)$ for each $K\in H^\star(G)$
and $\trans{c}(P_y)=a'(y)$ for each $y\in S$.
\item[(ii)] For each $y,y'\in S$ and $K\in H^\star(G)$,
$$a(K)+a'(y')-a'(y)\le \transplus{f}(K+P_{y'}-P_y+B^\star(G)).$$
\item[(iii)] For each $y,y'\in S$, $K\in H^\star(G)$ and a coboundary $R$,
if the copath $D=K+P_{y'}-P_y+R$ is simple, then
$$a(K)+a'(y')-a'(y)\le \transplus{f}(D).$$
\end{itemize}
Moreover, there is an algorithm that in time $O(n^2+(g+1)|S|n)$ either finds an $f$-circulation $c$ as described in (i),
or a simple copath $D=K+P_{y'}-P_y+R$ for some $y,y'\in S$, a coboundary $R$ and cocycle $K\in H^\star(G)$
such that $a(K)+a'(y')-a'(y)>\transplus{f}(D)$.
\end{theorem}
\begin{proof}
By Lemma~\ref{lemma-nomore}, (i) implies (ii), since $\trans{c}(K+P_{y'}-P_y)=a(K)+a'(y')-a'(y)$.
Moreover, (ii) clearly implies (iii), since $\transplus{f}(D)\ge \transplus{f}(K+P_{y'}-P_y+B^\star(G))$.
Hence, it remains to show that (iii) implies (i).  We do so by describing an algorithm
that either produces an $f$-circulation satisfying (i), or a simple copath $D$ showing that (iii) is false.

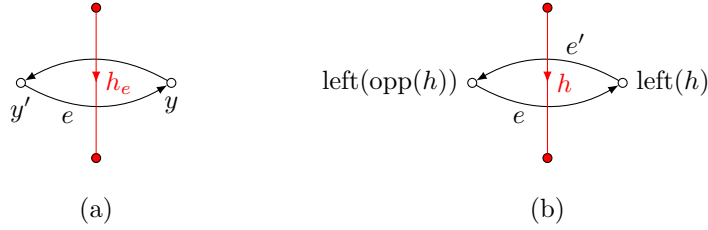
\begin{figure}
\begin{center}
\begin{tikzpicture}
\draw
(2,0) node[label=below:$y$,vtxw](y){}
(0,0) node[label=below:$y'$,vtxw](y'){}
(1,1) node[label=right:$ $,vtxr](x1){}
(1,-1) node[label=right:$ $,vtxr](x2){}
;
\draw[-latex](y') to[bend right] node[pos=0.3,below]{$e$} (y);
\draw[-latex](y) to[bend right] (y');
\draw[rededge] (x1)-- node[pos=0.5,right]{$h_e$} (x2) ;
\draw (1,-1.4) node[below]{(a)};

\begin{scope}[xshift = 6cm]
\draw (1,-1.4) node[below]{(b)};
\draw
(2,0) node[label=right:{left($h$)},vtxw](y){}
(0,0) node[label=left:{left(opp($h$))},vtxw](y'){}
(1,-1) node[label=right:$ $,vtxr](x1){}
(1,1) node[label=right:$ $,vtxr](x2){}
;
\draw[-latex](y') to[bend right] node[pos=0.3,below]{$e$} (y);
\draw[-latex](y) to[bend right] node[pos=0.3,above]{$e'$} (y');
\draw[rededge] (x2)-- node[pos=0.5,right]{$h$} (x1) ;
\end{scope}
\end{tikzpicture}
\end{center}
\caption{Notation in Theorem~\ref{thm-extcirc}.}\label{fig-ehe}
\end{figure}

By Lemma~\ref{lemma-extflow}, there always exists a $1$-cycle $b$ satisfying $\trans{b}(K)=a(K)$ for each $K\in H^\star(G)$
and $\trans{b}(P_y)=a'(y)$ for each $y\in S$; however, $b$ is not necessarily an $f$-circulation.
To obtain an $f$-circulation $c$, we modify $b$ by adding $\partial_2 L$ for a suitably chosen 2-chain $L$
such that $L[y]=0$ for every $y\in S$; by Observation~\ref{obs-transf}, this ensures that $\trans{c}(P)=\trans{b}(P)$
for every $P\in H^\star(G)\cup \{P_y:y\in S\}$.  The $2$-chain $L$ is chosen so that for every $y\in F(G)$,
$L[y]$ is the distance from $S$ to $y$ in a directed graph derived from the dual of $G$,
where the lengths of edges are carefully chosen to ensure that $0\le (b+\partial_2 L)[h]\le f[h]$ for each $h\in \hes(G)$
such that $f[h]\ge 0$.  Let us now work out the details of this idea.

Let us view the dual $G^\star$ of $G$ as a symmetrically oriented graph, and let $\ell$ be a function assigning
lengths to its directed edges as follows: For any edge $e=(y',y)\in E(G^\star)$,
let $h_e$ be the half-edge of $G$ dual to $e$ such that $\lft(h_e)=y$, see Figure~\ref{fig-ehe}(a).  Let
\begin{align}
\ell(e)=\begin{cases}
f[h_e]-b[h_e]&\text{ if $f[h_e]>0$}\\
-b[h_e]&\text{ if $f[h_e]\le 0$.}
\end{cases}
\label{eq:l}
\end{align}
\begin{figure}
\begin{center}
\includegraphics[page=15,scale=0.7]{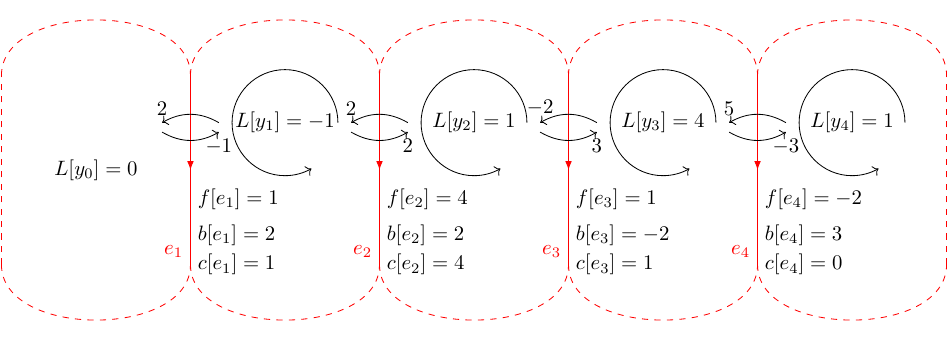}
\end{center}
\caption{Situation from Theorem~\ref{thm-extcirc}. It depicts 5 faces $y_0,\ldots,y_4$, where $y_0 \in S$ and $e_i$ is an edge such that
$\lft(e_i)=y_i$ and $\lft(\opp(e_i)) = y_{i-1}$, and the values of $f$, $b$, $c$ for each of the corresponding half-edges. The black edges crossing $e_i$ are edges of the
dual and the lengths $\ell$ of these edges as in \eqref{eq:l} are shown. The value of $L[y_i]$ is depicted in each face and
the directed arc shows the direction of $\partial_2 y_i$.}\label{fig-fbc}
\end{figure}
Let us remark that $\ell$ is (in general) neither symmetric nor antisymmetric, see Figure~\ref{fig-fbc}.

Let us first consider the case that $G^\star$ with this assignment of lengths contains a directed cycle $C$ of negative length,
and thus $\ell$ cannot be used to define a metric.  Let $D=\sum_{e\in E(C)} h_e$ and note that $D$ is a simple cocycle.  We have
\begin{align*}
\transplus{f}(D)-\trans{b}(D)&=\sum_{e\in E(C):f[h_e]>0} f[h_e] - \sum_{e\in E(C)} b[h_e]\\
&=\sum_{e\in E(C):f[h_e]>0} (f[h_e]-b[h_e]) + \sum_{e\in E(C):f[h_e]\le 0} (-b[h_e])\\
&=\sum_{e\in E(C)} \ell(e)<0.
\end{align*}
Moreover, since $D$ is a cocycle, we have $D=K+R$ for some $K\in H^\star(G)$ and a coboundary $R$, and
using Observation~\ref{obs-transv},
$$\transplus{f}(D)<\trans{b}(D)=\trans{b}(K)=a(K).$$
This shows that (iii) is false, with $y=y'=x$.

Hence, we can assume that there is no cycle of negative length in $G^\star$.  Suppose now that
there exists a directed path $P$ from some $y\in S$ to some $y'\in S$ in $G^\star$ of negative length,
and let $D=\sum_{e\in E(P)} h_e$.  Then $D$ is a simple copath from $y$ to $y'$, and we have
\begin{align*}
\transplus{f}(D)-\trans{b}(D)&=\sum_{e\in E(P):f[h_e]>0} (f[h_e]-b[h_e]) + \sum_{e\in E(P):f[h_e]\le 0} (-b[h_e])\\
&=\sum_{e\in E(P)} \ell(e)<0.
\end{align*}
Note that $D-P_{y'}+P_y$ is a cocycle, and thus $D-P_{y'}+P_y=R+K$ for some $K\in H^\star(G)$ and a coboundary $R$.
Hence, $D=K+P_{y'}-P_y+R$ and using Observation~\ref{obs-transv},
$$\transplus{f}(D)<\trans{b}(D)=\trans{b}(K)+\trans{b}(P_{y'})-\trans{b}(P_y)=a(K)+a'(y')-a'(y).$$
This again shows that (iii) is false.

Therefore, we can assume that $G^\star$ does not contain any path of negative length between elements of $S$.
Let $L$ be the $2$-chain such that for each $y\in F(G)$, $L[y]$ is the distance from $S$ to $y$
in $G^\star$ with edge lengths $\ell$; the distance is defined since there are no negative length cycles,
and we have $L[y]=0$ for every $y\in S$ since there are no negative length paths between elements of $S$.
Let $c=b+\partial_2 L$, see Figure~\ref{fig-fbc}.  We want to show $c$ is the desired $f$-circulation for (i). Since $\partial_2 L$ is a $1$-boundary and $b$ is a 1-cycle,
$c$ is a $1$-cycle.  By Observation~\ref{obs-transf}, we have $\trans{\partial_2 L}(K)=0$ for any cocycle $K$,
and thus $\trans{c}(K)=\trans{b}(K)=a(K)$ for each $K\in H^\star(G)$.  For $y\in S$,
we have
\begin{align*}
\trans{\partial_2 L}(P_y)&=\sum_{y'\in F(G)} L[y']\trans{\partial_2 y'}(P_y)=\sum_{y'\in F(G)} L[y'](\partial^\star_1 P_y)[y']\\
&=L[y]-L[x]=0,
\end{align*}
since $x,y\in S$.  Consequently, $\trans{c}(P_y)=\trans{b}(P_y)=a'(y)$.

Let us now consider a half-edge $h\in \hes(G)$, where $f[h]\ge 0$ without loss of generality,
let $e$ be the corresponding edge of $G^\star$ directed towards $\lft(h)$,
and let $e'$ be the edge opposite to $e$, see Figure~\ref{fig-ehe}(b).
We have
$$c[h]=b[h]+(\partial_2 L)[h]=b[h]+(L[\lft(h)]-L[\lft(\opp(h))]).$$
Since $L[\lft(h)]$ and $L[\lft(\opp(h))]$ are distances from $S$ to $\lft(h)$ and to $\lft(\opp(h))$,
triangle inequality gives
\begin{align*}
L[\lft(h)]&\le L[\lft(\opp(h))]+\ell(e)=L[\lft(\opp(h))]+f[h]-b[h]
\end{align*}
and
\begin{align*}
L[\lft(\opp(h))]&\le L[\lft(h)]+\ell(e')=L[\lft(h)]-b[-h]=L[\lft(h)]+b[h].
\end{align*}
Therefore,
$$-b[h]\le L[\lft(h)]-L[\lft(\opp(h))]\le f[h]-b[h],$$
and thus $0\le c[h]\le f[h]$.  Hence, $c$ is an $f$-circulation, and (i) holds.

The single-source shortest paths (or a negative cycle) from $S$ can be found using
Bellman-Ford algorithm in time $O(|V(G^\star)|\cdot|E(G^\star)|)=O(n^2)$.
Adding the complexity of the algorithm from Lemma~\ref{lemma-extflow}, we conclude
that the described procedure can be implemented with time complexity $O(n^2+(g+1)|S|)n)$.
\end{proof}

Let us note the following flow decomposition corollary.
\begin{corollary}\label{cor-decomp}
Let $G$ be a graph with a 2-cell drawing in an orientable surface.
Let $c_0$ and $c$ be 1-cycles in $G$ and let $k$ be a positive integer.
If $c-kc_0\in B_1(G)$, then there exist $c$-circulations $c_1$, \ldots, $c_k$
such that $c=c_1+\ldots+c_k$ and $c_i-c_0\in B_1(G)$ for $i\in\{1,\ldots,k\}$.
\end{corollary}
\begin{proof}
We prove the claim by induction on $k$.  If $k=1$, then we can set $c_1=c$.  Hence, we can assume $k\ge 2$.
It suffices to show that a $c$-circulation $c_k$ such that $c_k-c_0\in B_1(G)$ exists, as the
claim then follows by the induction hypothesis for $k-1$ and the 1-cycle $c-c_k$.

Let a homomorphism $a:H^\star(G)\to\mathbb{Z}$ be defined by setting $a(K)=\trans{c_0}(K)$ for each $K\in H^\star(G)$.
Note that $\trans{c}(K)=\trans{kc_0}(K)=ka(K)$ by Observation~\ref{obs-transf}.
Let $S=\{x\}$ for an arbitrary face $x$ of $G$, and let $a':S\to\mathbb{Z}$ be defined as $a'(x)=0$.
Since $c$ is a $c$-circulation, Lemma~\ref{lemma-nomore} implies 
that for any $K\in H^\star(G)$, we have
$$ka(K)=\trans{c}(K)\le \transplus{c}(K+B^\star(G)).$$
Since the right-hand side is non-negative and $k\ge 1$, this implies
$$a(K)\le \transplus{c}(K+B^\star(G)).$$
Hence, Theorem~\ref{thm-extcirc} implies the existence of a $c$-circulation $c_k$
such that $\trans{c_k}(K)=a(K)=\trans{c_0}(K)$ for each $K\in H^\star(G)$.
Therefore, $\trans{c_k-c_0}(K)=0$ for each $K\in H^\star(G)$, and thus
$c_k-c_0\in B_1(G)$ by Corollary~\ref{cor-testb}.
\end{proof}

We will actually need the dual form of this corollary.
By a \emph{$K$-cocirculation}, we mean a cocycle $K'\preceq K$.
\begin{corollary}\label{cor-decompdual}
Let $G$ be a graph with a 2-cell drawing in an orientable surface.
Let $K_0$ and $K$ be cocycles in $G$ and let $k$ be a positive integer.
If $K-kK_0\in B^\star(G)$, then there exist $K$-cocirculations $K_1$, \ldots, $K_k$
such that $K=K_1+\ldots+K_k$ and $K_i-K_0\in B^\star(G)$ for $i\in\{1,\ldots,k\}$.
\end{corollary}

Let us now make some observations that are useful when dealing with the condition (ii) from Theorem~\ref{thm-extcirc}.
\begin{observation}\label{obs-subaddit}
Let $G$ be a graph with a 2-cell drawing in an orientable surface
and let $f$ be a 1-chain in $G$.   All $1$-chains $K$ and $K'$ in $G$ satisfy
$$\transplus{f}(K+K'+B^\star(G))\le \transplus{f}(K+B^\star(G))+\transplus{f}(K'+B^\star(G)).$$
\end{observation}
\begin{proof}
Let $R,R'\in B^\star(G)$ be such that $\transplus{f}(K+B^\star(G))=\transplus{f}(K+R)$ and
$\transplus{f}(K'+B^\star(G))=\transplus{f}(K'+R')$.
Since $R+R'\in B^\star(G)$,
\begin{align*}
\transplus{f}(K+K'+B^\star(G))&\le\transplus{f}(K+K'+R+R')\\
&=\sum_{h\in\hes{G}:f[h]>0} f[h]\max(0,(K+K'+R+R')[h])\\
&\le\sum_{h\in\hes{G}:f[h]>0} f[h](\max(0,(K+R)[h])+\max(0,(K'+R')[h]))\\
&= \transplus{f}(K+R)+\transplus{f}(K'+R').
\end{align*}
\end{proof}
Moreover, Corollary~\ref{cor-decompdual} implies multiplicativity of $\transplus{f}$ over
cocycles.
\begin{lemma}\label{lemma-mult}
Let $G$ be a graph with a 2-cell drawing in an orientable surface
and let $f$ be a 1-chain in $G$.  For every positive integer $k$ and a cocycle $K$, we have
$$\transplus{f}(kK+B^\star(G))=k\cdot \transplus{f}(K+B^\star(G)).$$
\end{lemma}
\begin{proof}
The inequality $\transplus{f}(kK+B^\star(G))\le k\cdot \transplus{f}(K+B^\star(G))$ follows from Observation~\ref{obs-subaddit}.
For the converse inequality, consider a cocycle $K'$ such that $K'-kK\in B^\star(G)$ and
$\transplus{f}(K')=\transplus{f}(kK+B^\star(G))$.
By Corollary~\ref{cor-decompdual}, there exist $K'$-cocirculations $K_1$, \ldots, $K_k$
such that $K'=K_1+\ldots+K_k$ and $K_i-K\in B^\star(G)$ for $i\in\{1,\ldots,k\}$.
Note that since $K_i\preceq K'$, for every $h\in\hes(G)$ we have $K_i[h]\ge 0$ if $K'[h]>0$
and $K'[h]>0$ if $K_i[h]>0$,
and thus $\transplus{f}(K')=\sum_{i=1}^k \transplus{f}(K_i)$.  Therefore,
$$\transplus{f}(kK+B^\star(G))=\transplus{f}(K')=\sum_{i=1}^k \transplus{f}(K_i)\ge k\cdot \transplus{f}(K+B^\star(G)).$$
\end{proof}

\section{Polytopes of allowed homologies}\label{sec-allow}

The conditions (ii) and (iii) from Theorem~\ref{thm-extcirc} can be viewed as linear inequalities constraining
 the functions $a$ and $a'$, implying that the possible homologies of $f$-circulations form a
polytope in a space of bounded dimension.  In this section, we make this intuition precise.

\begin{figure}
\includegraphics[page=11,scale=1.1]{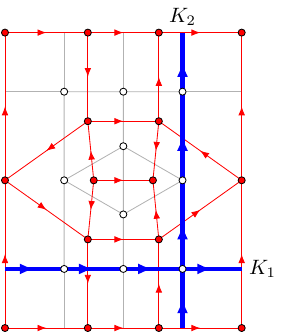}
\hskip 4em
\includegraphics[page=12,scale=0.7]{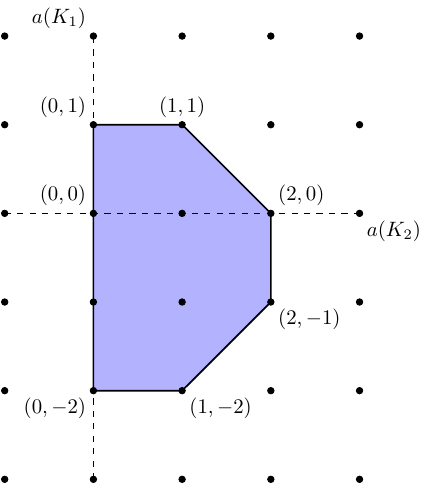}
\caption{A graph $G$ with a nowhere-zero flow $f$ drawn on the torus, a basis $Q=\{K_1,K_2\}$ of its cohomology group,
and the corresponding polytope $\PP_{G,F,Q}$ of allowed homologies of $f$-circulations.}\label{fig-polyex}
\end{figure}

\begin{figure}
\begin{center}
\includegraphics[page=9]{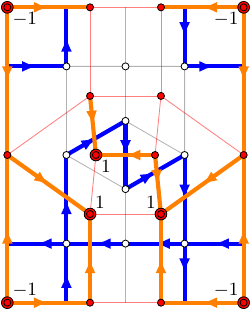}
\hskip 3em
\includegraphics[page=10]{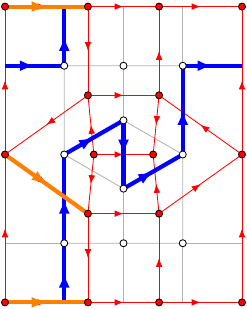}
\end{center}
\caption{Adding the coboundary $R \in B^\star(G)$ on the left to $K_1$ and $K_2$ gives the resulting cocycle on the right.
On the left the coefficients at vertices give the 0-chain $d$ such that $R=\partial^\star_2 d$, the orange edges show $R$
and the blue edges are the corresponding edges in $G^\star$.
On the right, the orange edges are the ones that are counted in $\transplus{f}( K_1 + K_2 + R)=2$. }
\label{fig:K1plusK2}
\end{figure}

Let $G$ be a graph with a 2-cell drawing in an orientable surface, let $Q$ be a basis of $H^\star(G)$,
and let $f$ be a $1$-chain in $G$.  Let $S$ be a non-empty subset of $F(G)$, let $x$ be an element of $S$ and
let $P$ be a function assigning to each $y\in S$ a copath $P(y)$ from $x$ to $y$, where $P(x)=0$.
For $a\in \mathbb{R}^Q$ and $z\in \mathbb{Z}^Q$, let $\langle z,Q\rangle=\sum_{K\in Q} z(K)K$ and $\langle z,a\rangle=\sum_{K\in Q} z(K)a(K)$.
Let us define
$$\PP_{G,f,Q,P}=\left\{\begin{array}{l}
(a,a')\in\mathbb{R}^{Q}\times \mathbb{R}^S:a'(x)=0,\\
\langle z,a\rangle+a'(y')-a'(y)\le \transplus{f}(\langle z,Q\rangle+P(y')-P(y)+B^\star(G))\\
\text{\hspace{38mm}for every $z\in\mathbb{Z}^Q$ and $y,y'\in S$}
\end{array}\right\};$$
and
$$\PP_{G,f,Q}=\{a\in\mathbb{R}^{Q}:\langle z,a\rangle\le \transplus{f}(\langle z,Q\rangle+B^\star(G))\text{ for every $z\in\mathbb{Z}^Q$}\}.$$
By Theorem~\ref{thm-extcirc}, there exists an $f$-circulation $c$ such that $\trans{c}(K)=a(K)$ for each $K\in Q$
and $\trans{c}(P(y))=a'(y)$ for each $y\in S$ if and only if $(a,a')\in\PP_{G,f,Q,P}$.
Figure~\ref{fig-polyex} shows an example of a graph on the torus with a nowhere-zero flow and the corresponding polytope $\PP_{G,f,\{K_1,K_2\}}$.
We also show how to obtain the constraint $a(K_1)+a(K_2) \leq 2$ for $\PP_{G,f,\{K_1,K_2\}}$. 
We start by fixing $z = (1,1)$. This gives constraint
\[
a(K_1)+a(K_2) \leq \transplus{f}( K_1 + K_2 + B^\star(G)).
\]
We find a particular $R \in B^\star(G)$ as depicted in Figure~\ref{fig:K1plusK2} on the left. 
This gives $\transplus{f}( K_1 + K_2 + B^\star(G))\le \transplus{f}( K_1 + K_2 + R)=2$.

Since $\PP_{G,f,Q,P}$ and $\PP_{G,f,Q}$ are given by infinite systems of inequalities, it is not
obvious that they are indeed polytopes (rather than more general convex shapes).  We show this together with other
useful properties in the following lemma.

\begin{lemma}\label{lemma-veras}
Let $G$ be a graph with a 2-cell drawing in an orientable surface, let $Q$ be a basis of $H^\star(G)$,
and let $f$ be a $1$-chain in $G$.  Let $S$ be a non-empty subset of $F(G)$, let $x$ be an element of $S$, and
let $P$ be a function assigning to each $y\in S$ a copath $P(y)$ from $x$ to $y$, where $P(x)=0$.
Then
\begin{itemize}
\item[(a)] $\PP_{G,f,Q,P}$ and $\PP_{G,f,Q}$ are polytopes determined by finite sets of inequalities,
\item[(b)] for each $K\in Q$ and $y\in S$, if $a\in \PP_{G,f,Q}$ or $(a,a')\in \PP_{G,f,Q,P}$, then
\begin{align*}
-\transplus{f}(-K)&\le a(K)\le \transplus{f}(K)\text{ and}\\
-\transplus{f}(-P(y))&\le a'(y)\le \transplus{f}(P(y)),
\end{align*}
\item[(c)] the vertices of $\PP_{G,f,Q}$ have integer coordinates, and
\item[(d)] for every $z\in\mathbb{Z}^Q$,
$$\transplus{f}(\langle z,Q\rangle+B^\star(G))=\max\{\langle z,a\rangle:a\in \PP_{G,f,Q}\}.$$
\end{itemize}
\end{lemma}
\begin{proof}
For any $y,y'\in S$ and any copath $D$ from $y$ to $y'$, note that $D-(P(y')-P(y))$ is a cocycle, and thus there
exists $z_D\in \mathbb{Z}^Q$ such that
$D-(P(y')-P(y))-\langle z_D,Q\rangle\in B^\star(G)$.   Similarly, for any cocycle $D$, let $z_D\in \mathbb{Z}^Q$
be such that $D-\langle z_D,Q\rangle\in B^\star(G)$.  Observe that
$$\PP=\{a\in\mathbb{R}^{Q}:\langle z_D,a\rangle\le \transplus{f}(D)\text{ for every simple cocycle $D$}\}$$
and
$$
\PP_P=\left\{\begin{array}{l}
(a,a')\in\mathbb{R}^{Q}\times \mathbb{R}^S:a'(x)=0,\\
\langle z_D,a\rangle+a'(y')-a'(y)\le \transplus{f}(D)\\
\text{\hspace{20mm}for every $y,y'\in S$ and a simple copath $D$ from $y$ to $y'$}
\end{array}\right\}$$
are polytopes determined by finite sets of inequalities, since there are only $3^{|E(G)|}$ simple $1$-chains in $G$.
Moreover, by the equivalence of (ii) and (iii) in Theorem~\ref{thm-extcirc},
we have $\PP_{G,f,Q}=\PP$ and $\PP_{G,f,Q,P}=\PP_P$.  Therefore (a) holds.

For each $K\in Q$, the definition of $\PP_{G,f,Q}$ gives
$$a(K)\le \transplus{f}(K+B^\star(G))\le \transplus{f}(K)$$
and
$$-a(K)=a(-K)\le \transplus{f}(-K+B^\star(G))\le \transplus{f}(-K)$$
for each $a\in \PP_{G,f,Q}$ or $(a,a')\in \PP_{G,f,Q,P}$ (by considering the case $y=y'=x$);
and for $(a,a')\in \PP_{G,f,Q,P}$, we similarly get (by considering $y$ together with $x$ and $K=0$)
$$a'(y)\le \transplus{f}(P(y)+B^\star(G))\le \transplus{f}(P(y))$$
and
$$-a'(y)\le \transplus{f}(-P(y)+B^\star(G))\le \transplus{f}(-P(y)).$$
Therefore (b) holds.

Consider now any $z\in\mathbb{Z}^Q$, and let $m=\max\{\langle z,a\rangle:a\in \PP_{G,f,Q}\}$.
Note that $m\le \transplus{f}(\langle z,Q\rangle+B^\star(G))$, since the definition of $\PP_{G,f,Q}$
includes the inequality $\langle z,a\rangle\le \transplus{f}(\langle z,Q\rangle+B^\star(G))$.
If $m<\transplus{f}(\langle z,Q\rangle+B^\star(G))$, then by Farkas lemma and (a),
the inequality $\langle z,a\rangle<\transplus{f}(\langle z,Q\rangle+B^\star(G))$ follows from a linear combination with non-negative rational coefficients of the constraints defining $\PP_{G,f,Q}$;
i.e., there exists a function
$\lambda:\mathbb{Z}^Q\to \mathbb{Z}_0^+$ with finite support and a positive integer $D$
such that 
$$\sum_{z'\in\supp(\lambda)} \frac{\lambda(z')}{D}z'=z$$
and
$$\sum_{z'\in\supp(\lambda)} \frac{\lambda(z')}{D}\transplus{f}(\langle z',Q\rangle+B^\star(G))<\transplus{f}(\langle z,Q\rangle+B^\star(G)).$$
However, by Lemma~\ref{lemma-mult} and Observation~\ref{obs-subaddit}, we have
\begin{align*}
D\cdot \transplus{f}(\langle z,Q\rangle+B^\star(G))&=\transplus{f}(\langle Dz,Q\rangle+B^\star(G))\\
&=\transplus{f}\left(\left\langle\sum_{z'\in\supp(\lambda)} \lambda(z')z',Q\right\rangle+B^\star(G)\right)\\
&\le \sum_{z'\in\supp(\lambda)} \lambda(z')\transplus{f}(\langle z',Q\rangle+B^\star(G)),
\end{align*}
which is a contradiction.  Therefore, $m=\transplus{f}(\langle z,Q\rangle+B^\star(G))$ and (d) holds.

Since $\max\{\langle z,a\rangle:a\in \PP_{G,f,Q}\}$ is an integer for every $z\in\mathbb{Z}^Q$,
it follows that all vertices of $\PP_{G,f,Q}$ have integer coordinates and (c) holds.
\end{proof}

Note that the minimal system of inequalities defining $\PP_{G,f,Q,P}$ may have exponential size,
and thus for algorithmic purposes, we cannot afford to represent the polytope explicitly.
A \emph{separation oracle} for a polytope $\PP\subseteq \mathbb{R}^d$ is an algorithm
that for an input $u\in \mathbb{Q}^d$ either decides that $u\in\PP$, or returns $z\in\mathbb{Q}^d$
such that the size of the binary encoding of $z$ is polynomial in the size of the binary encoding of $u$
and $\langle z,a\rangle < \langle z,u\rangle$ for every $a\in\PP$.
Theorem~\ref{thm-extcirc} can be used to obtain a separation oracle for $\PP_{G,f,Q,P}$.

\begin{lemma}\label{lemma-sepor}
Let $G$ be a graph of size $n$ with a 2-cell drawing in an orientable surface of Euler genus $g$, let $Q$ be a basis of $H^\star(G)$,
and let $f$ be a $1$-chain in $G$.  Let $S$ be a non-empty subset of $F(G)$, let $x$ be an element of $S$, and
let $P$ be a function assigning to each $y\in S$ a copath $P(y)$ from $x$ to $y$, where $P(x)=0$.
There exists a separation oracle for $\PP_{G,f,Q,P}$ with time complexity $O(n^2+(g+1)|S|n)$.
\end{lemma}
\begin{proof}
Consider the input $(u,u')\in \mathbb{Q}^Q\times\mathbb{Q}^S$.  If $u'(x)\neq 0$, we can return $(0,z')$ such that
$z'(x)=\text{sgn}(u'(x))$ and $z'(y)=0$ for all $y\in S\setminus \{x\}$.  Hence, suppose that $u(x)=0$.  We find a positive integer $\mu$
such that $(a_0,a'_0)=\mu\cdot(u,u')$ belongs to $\mathbb{Z}^Q\times\mathbb{Z}^S$.  Note that $a_0$ can be extended to a homomorphism from $H^\star(G)$
to $\mathbb{Z}$.  We apply the algorithm from Theorem~\ref{thm-extcirc} for $G$, paths $P_y=P(y)$ for $y\in S$, $a_0$ and $a'_0$,
and the $1$-chain $\mu f$.

If the outcome is a $\mu f$-circulation $c$ such that $\trans{c}(K)=a_0(K)$ for each $K\in H^\star(G)$
and $\trans{c}(P_y)=a_0'(y)$ for each $y\in S$, then by Lemma~\ref{lemma-nomore},
for every $z\in\mathbb{Z}^Q$ and $y,y'\in S$, we have
$$\langle z,a_0\rangle+a'_0(y')-a'_0(y)=\trans{c}(\langle z,Q\rangle+P(y')-P(y))\le \transplus{\mu f}(\langle z,Q\rangle+P(y')-P(y)+B^\star(G)),$$
and thus
$$\langle z,u\rangle+u'(y')-u'(y)\le \transplus{f}(\langle z,Q\rangle+P(y')-P(y)+B^\star(G))$$
and $(u,u')\in \PP_{G,f,Q,P}$.

Otherwise, the outcome is a simple copath $D=P(y')-P(y)+\langle z,Q\rangle+R$ for some $y,y'\in S$, $z\in\mathbb{Z}^Q$ and a coboundary $R$, satisfying
$$\transplus{\mu f}(\langle z,Q\rangle+P(y')-P(y)+B^\star(G))\le \transplus{\mu f}(D)<\langle z,a_0\rangle+a'_0(y')-a'_0(y).$$
Consequently, letting $z'(y')=1$, $z'(y)=-1$ (or $z'(y)=0$ if $y=y'$), and $z'(t)=0$ for $t\in S\setminus\{y,y'\}$,
for any $(a,a')\in \PP_{G,f,Q,P}$ we have
$$\langle (z,z'),(a,a')\rangle\le \transplus{f}(\langle z,Q\rangle+P(y')-P(y)+B^\star(G))<\langle (z,z'),(u,u')\rangle,$$
and thus we can return $(z,z')$.
\end{proof}

Let us remark that while we have included scaling by $\mu$ in the proof of Lemma~\ref{lemma-sepor}
to match the statement of Theorem~\ref{thm-extcirc}, it is easy to check that the algorithm from
Theorem~\ref{thm-extcirc} can be also directly used for (possibly non-integral) $u$ and $u'$.

Let $G$ be a graph with a 2-cell drawing in an orientable surface of Euler genus $g$, let $Q$ be a basis of $H^\star(G)$,
and let $f$ be a $1$-chain in $G$.  Let $S$ be a non-empty subset of $F(G)$, let $x$ be an element of $S$, and
let $P$ be a function assigning to each $y\in S$ a copath $P(y)$ from $x$ to $y$, where $P(x)=0$.
When we deal with the precoloring of an unbounded number of vertices, $|S|$ can be unbounded, and thus
the polytope $\PP_{G,f,Q,P}$ has unbounded dimension, preventing us from applying general integer programming
results.  To deal with this issue, for a fixed $a\in \mathbb{R}^Q$, let us define
$$\PP_{G,f,Q,a,P}=\left\{\begin{array}{l}
a'\in \mathbb{R}^S:a'(x)=0,\\
a'(y')-a'(y) \le \min_{z\in\mathbb{Z}^Q}\transplus{f}(\langle z,Q\rangle+P(y')-P(y)+B^\star(G))-\langle z,a\rangle\\
\text{\hspace{25mm}for every $y,y'\in S$}
\end{array}\right\}.$$
Clearly, we have $(a,a')\in\PP_{G,f,Q,P}$ if and only if $a\in \PP_{G,f,Q}$ and $a'\in \PP_{G,f,Q,a,P}$.
The dimension of $\PP_{G,f,Q}$ is $g$, and the polytope $\PP_{G,f,Q,a,P}$ is defined by polynomially many
inequalities that can be enumerated efficiently as shown in the following lemma.
\begin{lemma}\label{lemma-precas}
Let $\Sigma$ be an orientable surface of Euler genus $g$.  There exists an algorithm that,
given a graph $G$ of size $n$ with a 2-cell drawing in $\Sigma$, a basis $Q$ of $H^\star(G)$,
a $1$-chain $f$ in $G$, a non-empty subset $S$ of $F(G)$, an element $x$ of $S$,
a function $P$ assigning to each $y\in S$ a copath $P(y)$ from $x$ to $y$ (where $P(x)=0$),
and a vector $a\in \PP_{G,f,Q}\cap \mathbb{Z}^Q$, computes the right-hand sides of all inequalities
defining $\PP_{G,f,Q,a,P}$ in time $O(n^2+|S|n\log n)$.
\end{lemma}
\begin{proof}
Consider any $y,y'\in S$; we need to determine the minimum $\beta_{y,y'}$ of
$$\transplus{f}(\langle z,Q\rangle+P(y')-P(y)+B^\star(G))-\langle z,a\rangle$$
over all $z\in\mathbb{Z}^Q$.  Applying the equivalence between (i) and (ii) in Theorem~\ref{thm-extcirc}
with $\{y,y'\}$ playing the role of $S$ and with the copath $P(y')-P(y)$ from $y$ to $y'$,
it follows that $\beta_{y,y'}$ is the largest value such that there exists an $f$-circulation $c$
such that $\trans{c}(\langle z,Q\rangle)=\langle z,a\rangle$ for each $z\in\mathbb{Z}^Q$ and $\trans{c}(P(y')-P(y))=\beta_{y,y'}$.

Let $b$ be the $1$-cycle satisfying
$\trans{b}(\langle z,Q\rangle)=\langle z,a\rangle$ for each $z\in\mathbb{Z}^Q$
obtained using Lemma~\ref{lemma-extflow} with $\{x\}$ playing the role of $S$.
Let us remark that for any integer $\beta$, if we define $b_\beta=b+\beta\partial_2 y'$, then
$\trans{b_\beta}(\langle z,Q\rangle)=\langle z,a\rangle$ for each $z\in\mathbb{Z}^Q$ and $\trans{b_\beta}(P(y')-P(y))=\trans{b}(P(y')-P(y))+\beta$.
Let $\ell_\beta$ be the assignment of the lengths to the edges of the symmetric orientation of the dual $G^\star$
as defined in the proof of Theorem~\ref{thm-extcirc} with $b_\beta$ playing the role of $b$.  Note that for
any edge $e$ of $G^\star$ entering $y'$, letting $h_e$ be the corresponding dual edge such that $\lft(h_e)=y'$,
we have
$$\ell_\beta(e)=\begin{cases}
f[h_e]-b[h_e]-\beta&\text{ if $f[h_e]>0$}\\
-b[h_e]-\beta&\text{ if $f[h_e]\le 0$.}
\end{cases}$$
Since $a\in \PP_{G,f,Q}$, $G^\star$ with length assignment $\ell_\beta$ does not contain any negative cycle, and
$\beta_{y,y'}$ is the largest value of $\trans{b}(P(y')-P(y))+\beta$ such that there does not exist a path of negative length from $y'$ to $y$.
Denoting by $d_b(y,y')$ the distance from $y$ to $y'$ according to $\ell_0$, observe that the distance from $y$ to $y'$
according to $\ell_\beta$ is $d_b(y,y')-\beta$, and thus
$$\beta_{y,y'}=\trans{b}(P(y'))-\trans{b}(P(y))+d_b(y,y').$$
Hence, to determine the right-hand sides $\beta_{y,y'}$ for all $y,y\in S$, it suffices to determine the pairwise distances between all elements of $S$
in the length assignment $\ell_0$ for $\beta=0$.
Using Bellman-Ford algorithm followed by Johnson's reweighting and $|S|$ repetitions
of Dijkstra's algorithm, this can be done in time $O(n^2+|S|n\log n)$.  We can compute $b$ using Lemma~\ref{lemma-extflow}
and $\trans{b}(P(y))$ for each $y\in S$ in total time $O(n^2)$, and thus we obtain an algorithm with the desired time complexity.
\end{proof}

Finally, we will need a better understanding of the polytope $\PP_{G,f,Q}$ in the case that $f$ is a nowhere-zero flow.
For a $0$-chain $d$, a cocycle $K$, and a function $t:V(G)\to\mathbb{Z}$, let us define
\begin{align*}
\langle t,d\rangle&=\sum_{v\in V(G)} t(v)d(v)\\
\langle t,\partial^\star_2\rangle&=\sum_{v\in V(G)} t(v)\partial^\star_2 v\\
\transplus{d}(K+B^\star(G))&=\frac{1}{2}\min_{z'\in\mathbb{Z}^{V(G)}} |K + \langle z',\partial^\star_2\rangle|+\langle z',d\rangle\\
\PP_{G,d,Q}&=\{a\in\mathbb{R}^Q: \langle z,a\rangle\le\transplus{d}(\langle z,Q\rangle+B^\star(G))\text{ for every $z\in\mathbb{Z}^Q$}\}.
\end{align*}
The polytope $\PP_{G,d,Q}$ is an analogue of $\PP_{G,f,Q}$ defined in terms of the boundary,
and as we show next, $\PP_{G,f,Q}$ is a translation of $\PP_{G,\partial_1 f,Q}$, i.e.,
the polytope $\PP_{G,f,Q}$ is essentially the same for all nowhere-zero flows with the same boundary.

We are also going to need a result on the \emph{width} of these polytopes, defined as follows.
For a compact set $\PP\subseteq \mathbb{R}^d$ and a vector $z\in \mathbb{Z}^d$, let
$$w(\PP,z)=\max\{\langle z,a\rangle:a\in \PP\}-\min\{\langle z,a\rangle:a\in \PP\}.$$
Let $w(\PP)=\min\{w(\PP,z):z\in \mathbb{Z}^d\setminus\{0\}\}$.

\begin{lemma}\label{lemma-wide}
Let $G$ be a graph with a 2-cell drawing in an orientable surface,
let $Q$ be a basis of $H^\star(G)$, let $f$ be a nowhere-zero flow in $G$, and let $d=\partial_1 f$.
Then $\PP_{G,f,Q}$ is a translation of $\PP_{G,d,Q}$.
Moreover, there exist a cocycle $K\in Z^\star(G)\setminus B^\star(G)$ and a coboundary $R$ such that
$w(\PP_{G,f,Q})=(|K|+|K+R|-\trans{f}(R))/2$.
\end{lemma}
\begin{proof}
Note that for any cocycle $K_0$, Observation~\ref{obs-rev} gives
$$\transplus{f}(K_0)=(|K_0|+\trans{f}(K_0))/2.$$
Moreover, by Observation~\ref{obs-transv}, $\trans{f}(\langle z',\partial^\star_2\rangle)=\langle z',d\rangle$
for any $z'\in \mathbb{Z}^{V(G)}$.  Therefore,
\begin{align}
\transplus{f}(K_0+B^\star(G))&=\min_{z'\in \mathbb{Z}^{V(G)}} \transplus{f}(K_0+\langle z',\partial^\star_2\rangle)\nonumber\\
&=\frac{1}{2}\min_{z'\in \mathbb{Z}^{V(G)}} |K_0+\langle z',\partial^\star_2\rangle|+\trans{f}(K_0+\langle z',\partial^\star_2\rangle)\nonumber\\
&=\frac{1}{2}\trans{f}(K_0)+\frac{1}{2}\min_{z'\in \mathbb{Z}^{V(G)}} |K_0+\langle z',\partial^\star_2\rangle|+\langle z',d\rangle\nonumber\\
&=\frac{1}{2}\trans{f}(K_0)+\transplus{d}(K_0+B^\star(G)).\label{eq-tplus}
\end{align}
For any $a\in\mathbb{R}^Q$, define $a'\in\mathbb{R}^Q$ by letting $a'(K)=a(K)-\tfrac{1}{2}\trans{f}(K)$ for each $K\in Q$.
For any $z\in\mathbb{Z}^Q$, letting $K_0=\langle z,Q\rangle$ in (\ref{eq-tplus}) gives
\begin{align*}
\transplus{f}(\langle z,Q\rangle+B^\star(G))-\langle z,a\rangle&=\frac{1}{2}\trans{f}{\langle z,Q\rangle}+\transplus{d}(\langle z,Q\rangle+B^\star(G))-\langle z,a\rangle\\
&=\transplus{d}(\langle z,Q\rangle+B^\star(G))-\langle z,a'\rangle.
\end{align*}
Consequently, $a\in\PP_{G,f,Q}$ if and only if $a'\in \PP_{G,d,Q}$, and thus $\PP_{G,f,Q}$ is a translation of $\PP_{G,d,Q}$.

Let $z\in \mathbb{Z}^Q\setminus\{0\}$ be such that $w(\PP_{G,f,Q})=w(\PP_{G,f,Q},z)$, let $K_0=\langle z,Q\rangle$,
and let $R_1,R_2\in B^\star(G)$ be such that $\transplus{f}(K_0+B^\star(G))=\transplus{f}(K_0+R_1)$
and $\transplus{-f}(K_0+B^\star(G))=\transplus{-f}(K_0+R_2)$.
By Lemma~\ref{lemma-veras}(d) and Observation~\ref{obs-rev}, we have
\begin{align*}
\max\{\langle z,a\rangle:a\in \PP_{G,f,Q}\}&=\transplus{f}(K_0+B^\star(G))=\transplus{f}(K_0+R_1)\\
&=(|K_0+R_1|+\trans{f}(K_0+R_1))/2,
\end{align*}
\begin{align*}
\min\{\langle z,a\rangle:a\in \PP_{G,f,Q}\}&=-\max\{\langle -z,a\rangle:a\in \PP_{G,f,Q}\}\\
&=-\transplus{f}(-K_0+B^\star(G))=-\transplus{-f}(K_0+B^\star(G))\\
&=-\transplus{-f}(K_0+R_2)=-(|K_0+R_2|-\trans{f}(K_0+R_2))/2,
\end{align*}
and
\begin{align*}
w(\PP_{G,f,Q})&=(|K_0+R_1|+\trans{f}(K_0+R_1)+|K_0+R_2|-\trans{f}(K_0+R_2))/2\\
&=(|K_0+R_1|+|K_0+R_2|+\trans{f}(R_1-R_2))/2.
\end{align*}
Hence, for the cocycle $K=K_0+R_1$ and the coboundary $R=R_2-R_1$, we have
$w(\PP_{G,f,Q})=(|K|+|K+R|-\trans{f}(R))/2$.  Note that $K\not\in B^\star(G)$, since $z\neq 0$.
\end{proof}

\section{Circulations with prescribed modulo}\label{sec-nonem}

So far, we have considered the problem of determining whether there exists a circulation with a prescribed
value over given copaths and cocycles.  However, for the application in testing the existence of a homomorphism to an odd cycle, we need to know whether
there exists one with a given value modulo an odd integer.  The following lemma is used to deal with the precolored vertices
via the polytope $\PP_{G,f,Q,a,P}$.

\begin{lemma}\label{lemma-solvele}
There exists an algorithm that, given a finite set $S$, an element $x\in S$, a positive integer $m$,
a function $d:S^2\to \mathbb{Z}$, and a function $r:S\to \{0,1,\ldots,m-1\}$ such that $r(x)=0$, decides in time
$O(|S|^3)$ whether there exists a function $\ell:S\to\mathbb{Z}$ such that
\begin{itemize}
\item $\ell(x)=0$,
\item for every $y,y'\in S$, $\ell(y')-\ell(y)\le d(y,y')$, and
\item for every $y\in S$, $\ell(y)\equiv r(y)\pmod m$.
\end{itemize}
If such a function exists, the algorithm returns one.
\end{lemma}
\begin{proof}
For every $y,y'\in S$, let $d'(y,y')$ be the largest integer such that $d'(y,y')\le d(y,y')$
and $d'(y,y')\equiv r(y')-r(y)\pmod m$.  Clearly, if $a\equiv r(y)\pmod m$ and $b\equiv r(y')\pmod m$,
then $b-a\le d(y,y')$ if and only if $b-a\le d'(y,y')$.  Hence, it suffices to solve the problem with $d$ replaced by $d'$.
If $d'(y,y)<0$ for any $y\in S$, then $\ell$ does not exist; hence, suppose that $d'(y,y)\ge 0$ for every $y\in S$.

Consider the complete symmetrically oriented graph $K$ with vertex set $S$ and with each edge $(y,y')\in E(K)$
having length $d'(y,y')$.  If $K$ contains a cycle $y_0y_1y_2\ldots y_k$ of negative length, where $y_0=y_k$, then
the function $\ell$ does not exist, since otherwise we would have
$$0=\sum_{i=1}^k (\ell(y_i)-\ell(y_{i-1}))\le \sum_{i=1}^k d'(y_{i-1},y_i)<0.$$
Otherwise, let $\ell(y)$ be the distance from $x$ to $y$ in $K$.
By the triangle inequality, this ensures that $\ell(y')-\ell(y)\le d'(y,y')$ for every $y,y'\in S$.
Moreover, if $x=x_0,x_1,\ldots,x_t=y$ is a shortest path from $x$ to $y$ in $K$,
then
$$\ell(y)=\sum_{i=1}^t d'(x_{i-1},x_i)\equiv \sum_{i=1}^t (r(x_i)-r(x_{i-1}))=r(y)-r(x)=r(y)\pmod m.$$
Therefore, one can find $\ell$ or decide that it does not exist by applying
Bellman-Ford algorithm to $K$, with time complexity $O(|S|^3)$.
\end{proof}

Also, let us note the following simple observation, relevant because of Lemma~\ref{lemma-veras}.
\begin{observation}\label{obs-polymod}
For a polytope $\PP\subseteq \mathbb{R}^d$, a positive integer $m$, and a vector $r\in\mathbb\{0,\ldots,m-1\}^d$,
there exists a point $p\in \mathbb{Z}^d\cap \PP$ such that $p(i)\equiv r(i)\pmod m$ for $1\le i\le n$ if
and only if the polytope $(\PP-r)/m$ contains a point with integer coordinates.
\end{observation}

To take advantage of this observation, we use the following integer programming result of Dadush~\cite{Dadush12}.

\begin{theorem}[{Dadush~\cite[Algorithm 7.2]{Dadush12}}]\label{thm-ilp}
There exists a function $\gamma_0(d)=d^{(1+o(1))d}$ and an algorithm that, given
a convex set $\PP\subseteq [-n,n]^d$ described by a separation oracle with time complexity $T$, returns in time
$O(\gamma_0(d)T\pll n)$ a point in $\mathbb{Z}^d\cap \PP$ or decides no such point exists.
\end{theorem}

We can now combine the results from the previous sections to obtain the following key algorithm.
\begin{lemma}\label{lemma-maingen}
Let $\Sigma$ be an orientable surface of Euler genus $g$.
There exists a function $\gamma$ and an algorithm that, given a graph $G$ of size $n$ with a 2-cell drawing
on $\Sigma$, a basis $Q$ of $H^\star(G)$ formed by simple cocycles, a positive odd integer $m$, a non-empty subset $S$ of $F(G)$,
an element $x\in S$, a function $P$ assigning to each $y\in S$ a copath $P(y)$ from $x$ to $y$, where $P(x)=0$, 
a function $r:S\setminus\{x\}\to\{0,\ldots,m-1\}$, and a $0$-boundary $b$ divisible by $m$,
in time $O((|b|+1)n+\min(n^g(n^2+|S|^3),\gamma(|S|)n^2\pll n))$ either
\begin{itemize}
\item finds a nowhere-zero flow $f$ such that $\partial_1 f=b$, $m|\trans{f}(K)$ for every $K\in Q$, and 
$\trans{f}(P(y))\equiv r(y)\pmod m$ for each $y\in S\setminus\{x\}$, or
\item decides no such nowhere-zero flow exists.
\end{itemize}
\end{lemma}
\begin{proof}
Let us remark that the genus $g$ is a fixed constant and $|S|\le n$, and thus the time complexity of
the algorithms from Theorem~\ref{thm-extcirc} and Lemma~\ref{lemma-sepor} is $O(n^2)$.
Let $\gamma_0$ be the function from Theorem~\ref{thm-ilp}, and let us define $\gamma(s)=\gamma_0(g+s)$.

First, using Ford-Fulkerson algorithm and the algorithms from Lemmas~\ref{lemma-flow1} and~\ref{lemma-flow2}, we can
in time $O((|b|+1)n)$ find a nowhere-zero flow $f_0$ such that $\partial_1 f_0=b$, or decide no such flow exists.

Suppose we found such a flow $f_0$.  For $K\in Q$, let $r_0(K)=\frac{m+1}{2}\trans{f_0}(K)$,
for $y\in S\setminus\{x\}$ let $r_0'(y)=\frac{m+1}{2}(\trans{f_0}(P(y))-r(y))$, and let $r_0'(x)=0$.
By Corollary~\ref{cor-circ}, there exists a nowhere-zero flow $f$ with the properties described in
the statement of the lemma if and only if there exists an $f_0$-circulation $c$ such that
$\trans{c}(K)\equiv r_0(K)\pmod m$ for $K\in Q$ and $\trans{c}(P(y))\equiv r_0'(y)\pmod m$ for $y\in S$.

By Theorem~\ref{thm-extcirc} and Observation~\ref{obs-polymod}, such an $f_0$-circulation $c$ exists if and only if the polytope
$$\frac{\PP_{G,f_0,Q,P}-(r_0,r'_0)}{m}$$
contains a point with integer coordinates.  By Lemma~\ref{lemma-sepor}, Theorem~\ref{thm-ilp} and the bounds from Lemma~\ref{lemma-veras}(b),
we can decide whether such a point $(u,u')$ exists (and find it if this is the case)
in time $O(\gamma_0(g+|S|)n^2\pll n)=O(\gamma(|S|)n^2\pll n)$.
Given this point, we can then apply the algorithm from Theorem~\ref{thm-extcirc} with the homomorphism $a$ that
maps $Q$ to $mu+r_0$ and with $a'=mu'+r'_0$ to find the desired $f_0$-circulation~$c$.

Alternately, instead of using the integer programming, we can go over $O(n^g)$ points $u$ within the bounds given by Lemma~\ref{lemma-veras}(b)
such that $u-r_0$ is divisible by $m$, and for each of them use the algorithm from Theorem~\ref{thm-extcirc} to determine whether
$u$ belongs to $\PP_{G,f_0,Q}$.  For each such point $u$ belonging to $\PP_{G,f_0,Q}$, we can then find the right-hand sides of
the inequalities defining $\PP_{G,f,Q,u,P}$ in time $O(n^2+|S|n\log n)$ by Lemma~\ref{lemma-precas}, then use
Lemma~\ref{lemma-solvele} to find $u'\in \PP_{G,f,Q,u,P}$ such that $u'-r'_0$ is divisible by $m$ in
time $O(|S|^3)$.  The time complexity per each point is $O(n^2+|S|n\log n+|S|^3)=O(n^2+|S|^3)$.
If such a point $u'$ is found, we proceed with $(u,u')\in \PP_{G,f_0,Q,P}$ as described in the previous paragraph.
\end{proof}

This easily implies our main algorithmic result.
\begin{proof}[Proof of Theorem~\ref{thm-maingen}]
Let $m=|C|$ and let $G$ be the dual of $H$.  Without loss of generality, we can assume that $S$ is non-empty (otherwise add any face of $G$ to $S$ and
color it arbitrarily); let $x$ be any element of $S$.  For each $y\in S\setminus\{x\}$, let $r(y)=(\psi(y)-\psi(x))\bmod m$.
Let $Q$ be a basis of $H^\star(G)$ obtained using Observation~\ref{obs-genhom}.  By performing a depth-first search in
$H$ from $x$, we find a simple copath $P(y)$ from $x$ to $y$ for each $y\in S\setminus\{x\}$,
and we let $P(x)=0$. 

By Lemma~\ref{lemma-tutte} and Observation~\ref{obs-relev}, it suffices to go over all (at most $q^\star_C(H)$)
relevant $0$-boundaries $b$ divisible by $m$ and for each of them verify whether there exists a nowhere-zero flow $f$
with $\partial_1 f=b$, $m|\trans{f}(K)$ for every $K\in Q$, and $\trans{f}(P(y))\equiv r(y)\pmod m$ for every $y\in S\setminus\{x\}$;
this can be done using the algorithm from Lemma~\ref{lemma-maingen}.
\end{proof}

\section{Polytopes with no integer points}\label{sec-hollow}

In order to obtain sufficient conditions for 3-colorability, Observation~\ref{obs-polymod} suggests
that we need a sufficient condition ensuring that a polytope contains an integer
point.  A polytope that does not contain any integer points is called \emph{hollow}.
It is intuitively clear that a hollow polytope must be quite narrow to fit in between the integer points.
This was made precise by Kannan and Lov\'asz~\cite{Kannan88},
with an improved bound that we state below by Rudelson~\cite{rudelson},
and the best possible bound in the 2-dimensional case given by Hurkens~\cite{Hurkens1990BlowingUC}.

\begin{theorem}\label{thm-narrow}
For every positive integer $d$, there exists $\mu_d=O(d^{4/3})$ such that
if a bounded polytope $\PP\subseteq \mathbb{R}^d$ is hollow, then
$w(\PP)<\mu_d$.  Moreover, $\mu_2=1+2/\sqrt{3}$.
\end{theorem}

A set $\PP$ is \emph{centrally symmetric} if there exists a point $p$ such that $\PP=p-\PP$.
A better bound is known for centrally symmetric polytopes.

\begin{theorem}[Banaszczyk~\cite{Banaszczyk96}]\label{thm-narrow-censym}
For every positive integer $d$, there exists $\overline{\mu}_d=O(d\log d)$ such that
if a bounded centrally symmetric polytope $\PP\subseteq \mathbb{R}^d$ is hollow, then
$w(\PP)<\overline{\mu}_d$.
\end{theorem}

We are going to need the following standard observation about change of lattice bases.
\begin{lemma}\label{lemma-cbas}
Let $Z\subseteq \mathbb{R}^d$ be a compact set and let $A\in \mathbb{Z}^{d\times d}$ be a matrix with $|\det A|=1$.
Then $Z$ is hollow if and only if $A^TZ$ is, and $w(Z)=w(A^TZ)$.  Moreover, for any $c\in\mathbb{Z}^d$,
we have $w(A^TZ,A^{-1}c)=w(Z,c)$.
\end{lemma}
\begin{proof}
Since $|\det A|=1$, Cramer's rule implies that $A^{-1}\in \mathbb{Z}^{d\times d}$.  Hence, if $x\in A^TZ$ has integer coordinates,
the point $A^{-1T}x\in Z$ has integer coordinates as well.  Conversely, if $p\in Z$ has integer coordinates, then so does
the point $A^Tp\in A^TZ$.  Hence, $Z$ is hollow if and only if $A^TZ$ is hollow.

For any $z\in\mathbb{R}^d$, we have
$\langle c,z\rangle=\langle A^{-1}c,A^Tz\rangle$, and thus $w(A^TZ,A^{-1}c)=w(Z,c)$.
Hence,
\begin{align*}
    w(Z)&=\min\{w(Z,c):c\in \mathbb{Z}^d\setminus\{0\}\}\\
    &=\min\{w(A^TZ,A^{-1}c):c\in \mathbb{Z}^d\setminus\{0\}\}\\
    &=\min\{w(A^TZ,c'):c'\in \mathbb{Z}^d\setminus\{0\}\}=w(A^TZ).
\end{align*}
\end{proof}

We say that a polytope $\PP$ is \emph{$\tfrac{1}{3}$-integral} if the coordinates of each vertex of $\PP$ are integer multiples of $1/3$.
Due to Observation~\ref{obs-polymod}, when considering 3-colorability, we are interested in the polytope
$(\PP_{G,f,Q}-r)/3$, which is $\tfrac{1}{3}$-integral
by Lemma~\ref{lemma-veras}(c).  For $\tfrac{1}{3}$-integral polytopes, the bound from the dimension 2 case of Theorem~\ref{thm-narrow}
can be improved; this is relevant when we consider graphs drawn on the torus (a surface of Euler genus 2).

\begin{lemma}\label{lemma-narrow}
If a bounded $\tfrac{1}{3}$-integral polytope $\PP\subseteq \mathbb{R}^2$ is hollow,
then $w(\PP)<2$.
\end{lemma}
\begin{proof}
Suppose for a contradiction that $w(\PP)\ge 2$.
By Theorem~\ref{thm-narrow}, there exists $c\in \mathbb{Z}^2\setminus\{(0,0)\}$ such that $w(\PP,c)<1+2/\sqrt{3} < 2.1548$,
and since $w(\PP,c)$ is a multiple of $1/3$ and $w(\PP)\ge 2$, we conclude that $w(\PP,c)=2$.

We are going to first transform $\PP$ into a $\tfrac{1}{3}$-integral hollow polytope $\PP'$ such that $w(\PP')=w(\PP',(1,0))=2$,
i.e., the $x$-coordinates of all points of $\PP$ are contained in an interval of length two.
Next, we further transform it into a hollow $\tfrac{1}{3}$-integral polytope $\PP''$ such that $w(\PP'')=2$ and 
$\PP\subseteq [0,2+2/3] \times [0,4+1/3]$.  There are only finitely many $\tfrac{1}{3}$-integral polytopes
contained in $[0,2+2/3] \times [0,4+1/3]$, and we show that all of them have width less than two by
computer-assisted enumeration, thus obtaining a contradiction.

Let $c=(m,n)^T$.  Note that $m$ and $n$ are co-prime, as otherwise we would have $w(\PP,c/\gcd(m,n))<2$.
Therefore, there exist integers $\alpha$ and $\beta$ such that $\alpha m+\beta n=1$.
Let $$A=\begin{pmatrix}
m&-\beta\\
n&\alpha
\end{pmatrix},$$ so that $\det A=1$ and $A^{-1}c=(1,0)^T$.  Let $\PP'=A^T\PP$.  By Lemma~\ref{lemma-cbas},
the polytope $\PP'$ is hollow and $w(\PP')=w(\PP',(1,0)^T)=2$.  Moreover,
$\PP'$ is clearly a $\tfrac{1}{3}$-integral polytope.

Let $l=(l_1,l_2)$ and $r=(r_1,r_2)$ be vertices of $\PP'$ with $l_1$ minimum and $r_1$ maximum;
$r_1-l_1=w(\PP',(1,0)^T)=2$.
Let $b$ be the integer such that $-1\le r_2-l_2+2b<1$, and let
$$B=\begin{pmatrix}
1&b\\
0&1
\end{pmatrix}.$$
Note that $\det B=1$ and  $B^{-1}(1,0)^T=(1,0)^T$, and thus letting $\PP''=B^T\PP'$, we have $w(\PP'',(1,0)^T)=w(\PP',(1,0)^T)=2$.
Let
\begin{align*}
l'&=(l'_1,l'_2)^T=B^T(l_1,l_2)^T=(l_1,l_2+bl_1)^T\text{ and}\\
r'&=(r'_1,r'_2)^T=B^T(r_1,r_2)^T=(r_1,r_2+br_1)^T,
\end{align*}
and note that $r'_2-l'_2=r_2-l_2+b(r_1-l_1)=r_2-l_2+2b$, and thus we have $|l'_2-r'_2|\le 1$.

Since $r_1-l_1=2$ and $l_1$ and $r_1$ are multiples of $1/3$, there exists an integer $x$ such that $l_1<x<r_1$,
$x-l_1\ge 2/3$, and $r_1-x\ge 2/3$.
Let $(x,y)$ be the point on the line between $l'$ and $r'$ whose first coordinate is $x$;
since $\PP''$ is hollow, $y$ is not an integer.
Let $L$ and $R$ be the open cones with apices $l'$ and $r'$, respectively, and rays passing through the points $(x,\lfloor y\rfloor)$ and
$(x,\lceil y\rceil)$.  Let $B$ be the band $\{(u_1,u_2):l_1\le u_1\le r_1,u_2\in\mathbb{R}\}$, see Figure~\ref{fig-LemNarrowClaim}.
We claim that $\PP''\subseteq (L\cup R)\cap B$.  Indeed, consider any point $u=(u_1,u_2)\in \PP''$.  Clearly $u\in B$.  Without loss of
generality, we can assume that $u_1\ge x$.  Let $(x,y')$ be the point on the line between $l'$ and $u$ whose first coordinate is $x$.
By convexity, $(x,y)$, $(x,y')$, and the whole segment between them lies in $\PP''$.  Since $\PP''$ is hollow, this segment does not
contain any point with integer coordinates, and thus $\lfloor y\rfloor < y'<\lceil y\rceil$, implying that $(x,y')$, and thus also $u$,
lies in the open cone $L$.

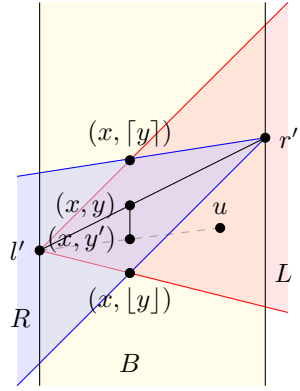
\begin{figure}
\center{
\begin{tikzpicture}[scale=.3]
\fill[yellow!10] (1,0) -- (11,0) -- (11,17) -- (1,17) -- (1,0);
\fill[red!20, opacity=.5] (1, 6) -- (12, 17) -- (12,3.25) -- (1,6);
\draw[red] (1,6) -- (12, 17);
\draw[red] (1,6) -- (12, 3.25);
\fill[blue!20, opacity=.5] (11,11) -- (0, 9.2857) -- (0,0) -- (11,11);
\draw[blue] (11,11) -- (0, 9.2857);
\draw[blue] (11,11) -- (0,0);
\draw(1,0) -- (1, 17);
\draw(11,0) -- (11, 17);
\draw(1,6) -- (11, 11);
\draw(5,8) -- (5, 6.5);
\draw[dashed, opacity=.3] (1,6) -- (9,7);
\node[vtx, label=left:{$l'$}] at (1,6) {};
\node[vtx, label=right:{$r'$}] at (11,11) {};
\node[vtx, label=left:{$(x,y)$}] at (5, 8) {};
\node[vtx, label=above:{$(x, \lceil y \rceil)$}] at (5,10) {};
\node[vtx, label=below:{$(x, \lfloor y \rfloor)$}] at (5,5) {};
\node[vtx, label=above:{$u$}] at (9, 7) {};
\node[vtx, label=left:{$(x,y')$}] at (5, 6.5) {};
\draw (5,1) node {$B$};
\draw (11,5) node[right] {$L$};
\draw (1,3) node[left] {$R$};
\end{tikzpicture}
}
\caption{Situation in Lemma~\ref{lemma-narrow}.}
\label{fig-LemNarrowClaim}
\end{figure}

Let $\delta=\lceil y\rceil-y$, and note that $y-\lfloor y\rfloor=1-\delta$.  Observe that
$\sup\{\langle (0,1)^T,z\rangle:z\in L\cap B\}$ is either $l'_2$ or
$$r'_2+\frac{r_1-l_1}{x-l_1}(\lceil y\rceil-y)=r'_2+\frac{2\delta}{x-l_1}.$$
Analogously, $\sup\{\langle (0,1)^T,z\rangle:z\in R\cap B\}$ is either $r'_2$ or $l'_2+\frac{2\delta}{r_1-x}$.
Therefore,
$$\max\{\langle (0,1)^T,z\rangle:z\in \PP''\}<\max\left(r'_2+\frac{2\delta}{x-l_1},l'_2+\frac{2\delta}{r_1-x}\right).$$
Analogously,
$$\min\{\langle (0,1)^T,z\rangle:z\in \PP''\}>\min\left(r'_2-\frac{2(1-\delta)}{x-l_1},l'_2-\frac{2(1-\delta)}{r_1-x}\right).$$
Consequently,
\begin{align*}
w(\PP'',(0,1)^T)<\max\Bigl(&\frac{2\delta}{x-l_1}+\frac{2(1-\delta)}{x-l_1},
\frac{2\delta}{r_1-x}+\frac{2(1-\delta)}{r_1-x},\\
&r'_2-l'_2+\frac{2\delta}{x-l_1}+\frac{2(1-\delta)}{r_1-x},\\
&l'_2-r'_2+\frac{2\delta}{r_1-x}+\frac{2(1-\delta)}{x-l_1}\Bigr).
\end{align*}
Since $x-l_1\ge 2/3$ and $r_1-x\ge 2/3$, the first two terms of the maximum are smaller than or equal to $3$.
By symmetry, we can assume that $x-l_1\le r_1-x$; then the third term is non-decreasing in $\delta$,
and thus
$$r'_2-l'_2+\frac{2\delta}{x-l_1}+\frac{2(1-\delta)}{r_1-x}<r'_2-l'_2+\frac{2}{x-l_1}\le 1+3=4.$$
Analogously, the fourth term is smaller than $4$.
Since $w(\PP'',(0,1)^T)$ is a multiple of $1/3$, we conclude that $w(\PP'',(0,1)^T)\le 4-1/3$.

We can shift $\PP''$ by an integer vector if necessary to ensure that $0\le u_1\le 2+2/3$ and $0\le u_2\le 4+1/3$ for every $(u_1,u_2)\in \PP''$.
However, there are only finitely many $\tfrac{1}{3}$-integral polytopes with this property, and by computer-assisted enumeration,
we verified that for all of them, either $w(\PP'')<2$ or $\PP''$ contains an integer point\footnote{We wrote two independent programs in SageMath and C++. The programs used for this
verification can be found as ancillary files at arXiv posting of this paper.}.  This is a contradiction.
\end{proof}

To apply these results, we need to estimate the width of the polytope $\PP_{G,f,Q}$ in terms of the edgewidth of $G^\star$.
First, let us give several auxiliary results.
Let $G$ be a graph with a 2-cell drawing in an orientable surface $\Sigma$.
For a closed directed walk $W$ in $G^\star$, let $\wtococ{W}$ denote the corresponding cocycle.
For a cocycle $K$ in $G$, the \emph{support} $\vec{K}$ of $K$ is the directed
graph with vertex set $\{\lft(h):h\in \hes(G),K[h]\neq 0\}$ and containing $K[h]$ edges
from $\lft(\opp(h))$ to $\lft(h)$ for each half-edge $h$ such that $K[h]>0$.
The \emph{undirected support} $\overline{K}$ of $K$ is the undirected graph with the same vertex set as $\vec{K}$
and with edges $\{\lft(\opp(h)),\lft(h)\}$ for each half-edge $h$ such that $K[h]>0$; hence,
$\overline{K}$ is a subgraph of $G^\star$.  We view both $\vec{K}$ and $\overline{K}$ as drawn in $\Sigma$,
with their drawing inherited from $G^\star$.  
A (directed) cycle drawn in $\Sigma$ is \emph{separating} if deleting it from $\Sigma$ disconnects the surface,
and \emph{non-separating} otherwise.
\begin{observation}\label{obs-cobound}
Let $G$ be a graph with a 2-cell drawing in an orientable surface.
A (directed) cycle $C$ in $G^\star$ is separating if and only if the corresponding cocycle $\wtococ{C}$ is a coboundary.
\end{observation}

\begin{observation}\label{obs-nonsep}
Let $G$ be a graph with a 2-cell drawing in an orientable surface.
If a cocycle $K$ is not a coboundary, then there exists a non-separating directed cycle $C$ in $G^\star$
such that $\wtococ{C}\preceq K$.
\end{observation}
\begin{proof}
Since $K$ is a cocycle, the directed graph $\vec{K}$ is Eulerian, and thus
it can be expressed as an edge-disjoint union of directed cycles $C_1$, \ldots, $C_m$.
Since $K=\wtococ{C_1}+\cdots+\wtococ{C_m}$ is not a coboundary, we can assume that $\wtococ{C_1}$
is not a coboundary, and by Observation~\ref{obs-cobound}, $C_1$ is a non-separating cycle.
We clearly have $\wtococ{C_1}\preceq K$.
\end{proof}

For a cocycle $R$ in a graph $G$ with a 2-cell drawing in an orientable surface of non-zero Euler genus $g$,
let $\nu(R)$ denote the size of a smallest subset $X$ of $E(G^\star)$
such that $\overline{R}+X$ contains a non-separating cycle as a subgraph.
By Observation~\ref{obs-genhom}, the group $H^\star(G)$ is isomorphic to $\mathbb{Z}^g$ and thus non-trivial.
Hence, there exists a cocycle in $G$ that is not a coboundary.
Observation~\ref{obs-nonsep} implies that $G^\star$ contains a non-separating cycle,
and consequently, $\nu(R)$ is finite.  For a $1$-chain $f$, let $\mu_f(R)=2\nu(R)+|R|-\trans{f}(R)$.
In the following lemma, we show that the bound on the width of $\PP_{G,f,Q}$ given in Lemma~\ref{lemma-wide} can be simplified in terms of this notion.
\begin{lemma}\label{lemma-optk}
Let $G$ be a graph with a 2-cell drawing in an orientable surface other than the sphere and let $R$ be a coboundary.
Then
$$\min\{|K|+|K+R|-|R|:K\in Z^\star(G)\setminus B^\star(G)\}=2\nu(R).$$
Hence, if $f$ is a nowhere-zero flow in $G$ and $Q$ is a basis of  $H^\star(G)$,
then there exists a coboundary $R$ such that $w(\PP_{G,f,Q})\ge \tfrac{1}{2}\mu_f(R)$.
\end{lemma}
\begin{proof}
Observe that any real numbers $a$ and $b$ satisfy $|a|+|a+b|-|b|\ge 0$, and if
$b=0$, then $|a|+|a+b|-|b|=2|a|$.  Hence,
\begin{align*}
|K|+|K+R|-|R|&=\sum_{h\in \ohes(G)} |K[h]|+|(K+R)[h]|-|R[h]|\\
&\ge 2\sum_{h\in \ohes(G), R[h]=0} |K[h]|\ge 0
\end{align*}
for any cocycle $K$.
Consider a cocycle $K\in Z^\star(G)\setminus B^\star(G)$, let $C'_K$ be a non-separating directed cycle
in $G^\star$ such that $\wtococ{C'_K}\preceq K$ which exists by Observation~\ref{obs-nonsep},
and let $C_K$ be the underlying undirected cycle of $C'_k$.  Then
$$|K|+|K+R|-|R|\ge 2\sum_{h\in \ohes(G), R[h]=0} |K[h]|\ge 2|E(C_K)\setminus E(\overline{R})|\ge 2\nu(R).$$
Therefore, $\min\{|K|+|K+R|-|R|:K\in Z^\star(G)\setminus B^\star(G)\}\ge 2\nu(R)$,
and it suffices to show that there exists $K\in Z^\star(G)\setminus B^\star(G)$
such that $|K|+|K+R|-|R|=2\nu(R)$.

If $\vec{R}$ contains a directed non-separating cycle $C$ as a subgraph, then let $K'=\wtococ{C}$.
Note that $K'$ is not a coboundary by Observation~\ref{obs-cobound}, and that $K'\preceq R$.
Letting $K=K'-R$, we have
$$|K|+|K+R|-|R|=|K'-R|+|K'|-|R|=|R-K'|+|K'|-|R|=0=2\nu(R).$$
Hence, we can assume that all directed cycles in $\vec{R}$ are separating.

Suppose now that $X$ is a smallest set of edges of $G^\star$ such that $\overline{R}+X$ contains a non-separating cycle
as a subgraph, and let $C$ be such a cycle oriented in one of the two possible directions (chosen arbitrarily).
We view the edges of $X$ as directed along $C$, and let $Y$ be the $1$-chain formed by the sum of the corresponding
half-edges of $G$.  By Observation~\ref{obs-cobound}, $Q=\wtococ{C}$ is not a coboundary.  Let $W$ be the closed walk obtained from $C$
as follows:  For each edge $e$ of $C-X$ that is directed oppositely to a corresponding edge $e'$ of $\vec{R}$,
let $C_e$ be a cycle in $\vec{R}$ containing $e'$ (which exists, since $\vec{R}$ is Eulerian), and replace $e$ by the walk $C_e-e'$.
Note that this corresponds to the addition of $\wtococ{C_e}$ to $Q$.  
Since all directed cycles in $\vec{R}$ are separating, Observation~\ref{obs-cobound} implies that $Q'=\wtococ{W}$ is obtained
from $Q$ by adding coboundaries, and thus $Q'$ is not a coboundary.  By Observation~\ref{obs-nonsep},
there exists a directed non-separating cycle $C'$ such that letting $K'=\wtococ{C'}$, we have $K'\preceq Q'$.
By the minimality of $X$, we have $X\subseteq E(C')$.
Note that $W-X$ is a union of walks in $\vec{R}$, and thus $R[h]>0$ for each half-edge $h\not\in X$ such that $Q'[h]>0$.
Since $C'$ is a cycle and $K'=\wtococ{C'}\preceq Q'$, we have $0\le K'[h]\le 1\le R[h]$ for each such edge $h$, and thus
$K'-Y\preceq R$.  Letting $K=K'-R$, we have
\begin{align*}
|K|+|K+R|-|R|&=|K'-R|+|K'|-|R|\\
&=(|Y|+|R-(K'-Y)|)+(|K'-Y|+|Y|)-|R|\\
&=2|Y|=2|X|=2\nu(R).
\end{align*}
We conclude that $\min\{|K|+|K+R|-|R|:K\in Z^\star(G)\setminus B^\star(G)\}=2\nu(R)$.

Consider now a nowhere-zero flow $f$ in $G$ and a basis $Q$ of $H^\star(G)$.
By Lemma~\ref{lemma-wide}, there exists a cocycle $K\in Z^\star(G)\setminus B^\star(G)$ and a coboundary $R$ such that
\begin{align*}
2w(\PP_{G,f,Q})&=|K|+|K+R|-\trans{f}(R)\\
&=(|K|+|K+R|-|R|)+|R|-\trans{f}(R)\\
&\ge 2\nu(R)+|R|-\trans{f}(R)=\mu_f(R).
\end{align*}
\end{proof}

We also need the following standard topological observation on non-separating cycles.

\begin{figure}
\centering{
\includegraphics[width=12cm]{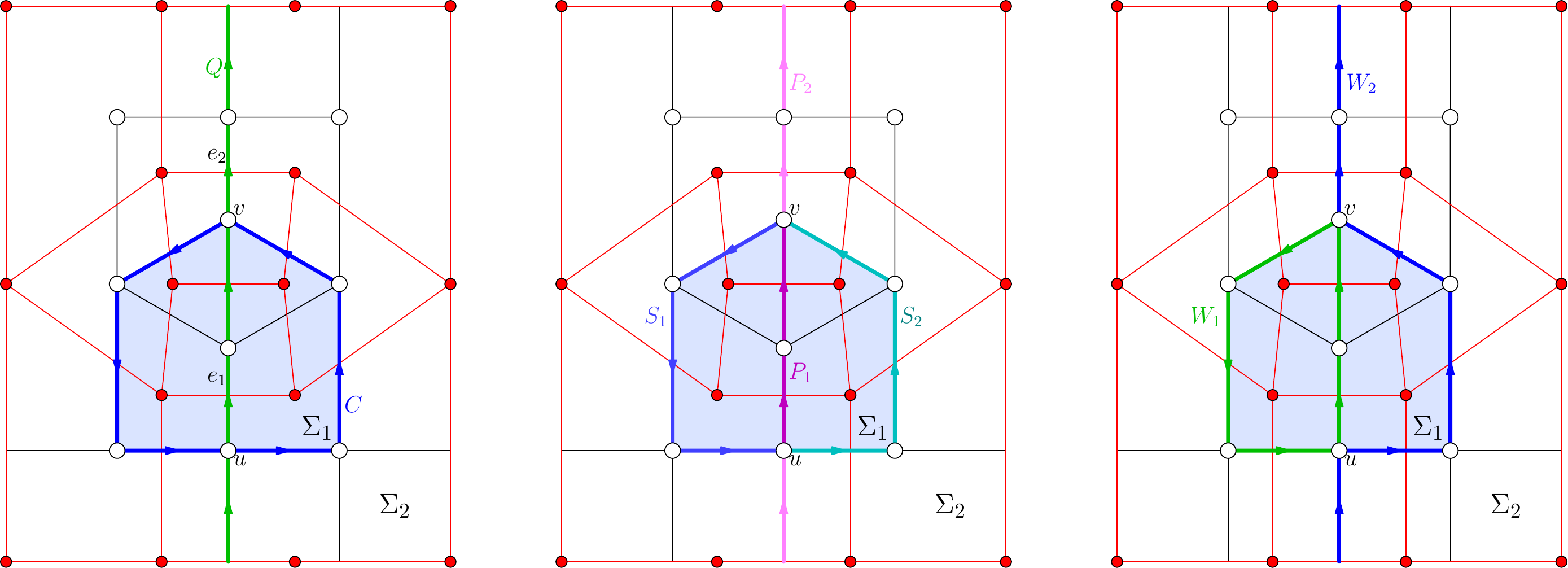}}
\caption{The situation in Lemma~\ref{lemma-minxsep}.}
\label{fig-LemMinxsep}
\end{figure}

\begin{lemma}\label{lemma-minxsep}
Let $G$ be a graph with a 2-cell drawing in an orientable surface $\Sigma$, let $C$ be a separating cycle
in $G^\star$, and let $\Sigma_1$ and $\Sigma_2$ be the connected parts of $\Sigma-C$.
For every non-separating cycle $Q$ in $G^\star$, there exists a non-separating cycle $Q'$
with $E(Q')\subseteq E(Q)\cup E(C)$ such that $Q'\subseteq \overline{\Sigma_1}$ or $Q'\subseteq \overline{\Sigma_2}$.
\end{lemma}
\begin{proof}
We prove the claim by the induction on $|E(Q)\setminus E(C)|$.
We can assume that there exist edges $e_1\in E(Q)\cap \Sigma_1$ and $e_2\in E(Q)\cap \Sigma_2$,
as otherwise we can set $Q'=Q$.  We orient the cycles $Q$ and $C$ in one of the two possible directions
arbitrarily.  There exist distinct vertices $u,v\in V(C)\cap V(Q)$ such that
$Q$ is the concatenation of a path $P_1$ from $u$ to $v$ and a path $P_2$ from $v$ to $u$,
where $e_1\in E(P_1)$ and $e_2\in E(P_2)$.
Moreover, $C$ is the concatenation of a path $S_1$ from $v$ to $u$ and $S_2$ from $u$ to $v$.
For $i\in\{1,2\}$, let $W_i$ be the closed walk obtained as the concatenation of $P_i$ and $S_i$, see Figure~\ref{fig-LemMinxsep}.
By Observation~\ref{obs-cobound}, $\wtococ{Q}$ is not a coboundary and $\wtococ{C}$ is, and
since $\wtococ{W_1}+\wtococ{W_2}=\wtococ{Q}+\wtococ{C}$, we can assume that $\wtococ{W_2}$ is not
a coboundary.  By Observation~\ref{obs-nonsep}, there exist a non-separating cycle $Q_2$ such that
$\wtococ{Q_2}\preceq \wtococ{W_2}$.  This implies that $E(Q_2)\subseteq E(P_2)\cup E(S_2)\subseteq (E(Q)\cup E(C))\setminus \{e_1\}$,
and thus $|E(Q_2)\setminus E(C)|<|E(Q)\setminus E(C)|$.
By the induction hypothesis, there exists
a non-separating cycle $Q'$ with $E(Q')\subseteq E(Q_2)\cup E(C)\subseteq E(Q)\cup E(C)$
such that $Q'\subseteq \overline{\Sigma_1}$ or $Q'\subseteq \overline{\Sigma_2}$,
as required.
\end{proof}

Let $R=\sum_{v\in V(G)} \alpha_v\partial^\star_2 v$ be a coboundary.  The \emph{span} of $R$ is
$\max\{\alpha_v:v\in V(G)\}-\min\{\alpha_v:v\in V(G)\}$.  Let us remark that the span of $R$
is independent of the way $R$ is expressed, since if $R$ can also be expressed as $\sum_{v\in V(G)} \alpha'_v\partial^\star_2 v$,
then there exists an integer $\delta$ such that $a'_v=a_v+\delta$ for every $v\in V(G)$.
\begin{lemma}\label{lemma-span}
Let $G$ be a graph with a 2-cell drawing in an orientable surface $\Sigma$ other than the sphere and let $R$ be a coboundary.
If $R$ has span at least three, then there exists a coboundary $R'\preceq R$ such that $R'\neq 0$
and $\nu(R-R')=\nu(R)$.
\end{lemma}
\begin{proof}
Let $R=\sum_{v\in V(G)} \alpha_v\partial^\star_2 v$.  Since $\sum_{v\in V(G)} \partial^\star_2 v=0$ and $R$ has span at least three,
we can without loss of generality assume that $\min\{\alpha_v:v\in V(G)\}\le -1$ and $\max\{\alpha_v:v\in V(G)\}\ge 2$.
For an integer $k\ge -1$, let $R_k=\sum_{v\in V(G):\alpha_v>k} \partial^\star_2 v$; note that
$R_k\preceq R$, since if $R_k[h]>0$ for a half-edge $h$, then $R_k[h]=1$, $\alpha_{\tgt(h)}>k$, $\alpha_{\tgt(\opp(h))}\le k$,
and $R[h]=\alpha_{\tgt(h)}-\alpha_{\tgt(\opp(h))}\ge 1$.
Moreover, $0\neq R_k\neq R$ for $k\in\{-1,0,1\}$, since $R$ has span at least three.
If $\overline{R_0}$ contains a non-separating cycle, then $\nu(R)=\nu(R_0)=0$, and thus we can set $R'=R-R_0$.
Hence, assume that all cycles in $\overline{R_0}$ are separating.

Let $X$ be a set of edges of $G^\star$ such that $|X|=\nu(R)$ and $\overline{R}+X$ contains
a non-separating cycle $Q$.  Choose such a cycle $Q$ with $|E(Q)\setminus E(\overline{R_0})|$
minimum.  Consider any cycle $C$ in $\overline{R_0}$.  We claim that all edges of $E(Q)\setminus E(\overline{R_0})$
are drawn in the same connected part of $\Sigma-C$; indeed, otherwise the cycle $Q'$ obtained by Lemma~\ref{lemma-minxsep}
would contradict the minimality of $|E(Q)\setminus E(\overline{R_0})|$.  Since this claim holds for all cycles in $\overline{R_0}$,
we conclude that there exists a face $x$ of $\overline{R_0}$ such that all edges of $E(Q)\setminus E(\overline{R_0})$
are drawn in $x$.  Note that either $\alpha_v\le 0$ for each $v\in V(G)$ drawn in $x$, or $\alpha_v>0$ for each $v\in V(G)$
drawn in $x$, since any vertices $v',v''\in V(G)$ such that $\alpha_{v'}\le 0$ and $\alpha_{v''}>0$ are necessarily
drawn in different faces of $\overline{R_0}$.  By symmetry, we can assume that the former is the case.
Consequently, $(E(Q)\setminus E(\overline{R_0}))\cap E(\overline{R_1})=\emptyset$,
and thus $\nu(R-R_1)=\nu(R)$.  Therefore, we can set $R'=R_1$.
\end{proof}

\begin{corollary}\label{cor-del}
Let $G$ be a graph with a 2-cell drawing in an orientable surface other than the sphere and let $f$ be a nowhere-zero flow in $G$.
If $R$ is a coboundary with $\mu_f(R)$ minimum and subject to that with $|R|$ minimum, then $R$ has span at most two.
\end{corollary}
\begin{proof}
Otherwise, let $R'\preceq R$ be the non-zero coboundary obtained by Lemma~\ref{lemma-span}.
We have
\begin{align*}
\mu_f(R-R')&=2\nu(R-R')+|R-R'|-\trans{f}(R-R')\\
&=2\nu(R)+|R|-|R'|-(\trans{f}(R)-\trans{f}(R'))\\
&=\mu_f(R)+\trans{f}(R')-|R'|\le \mu_f(R),
\end{align*}
and since $|R-R'|<|R|$, this contradicts the minimality of $R$.
\end{proof}

Finally, we are ready to bound the width of $\PP_{G,f,Q}$ in terms of the edgewidth of $G$.

\begin{corollary}\label{cor-width}
Let $G$ be a graph with a 2-cell drawing in an orientable surface other than the sphere,
let $Q$ be a basis of $H^\star(G)$, and let $f$ be a nowhere-zero flow in $G$.
Then $w(\PP_{G,f,Q})\ge \ew(G^\star)-|\partial_1 f|/2$.
\end{corollary}
\begin{proof}
By Lemma~\ref{lemma-optk}, there exists a coboundary $R$ such that $w(\PP_{G,f,Q})\ge \tfrac{1}{2}\mu_f(R)$.
Let $R_0$ be a coboundary with $\mu_f(R_0)$ minimum; by Corollary~\ref{cor-del}, we can assume that $R_0$ has
span at most two.  Hence, we can write $R_0=\sum_{v\in V(G)} \alpha_v\partial^\star_2 v$, where $-1\le \alpha_v\le 1$
for each $v\in V(G)$.  By Lemma~\ref{lemma-optk}, there exists a cocycle $K_0\in Z^\star(G)\setminus B^\star(G)$
such that $2\nu(R_0)=|K_0|+|K_0+R_0|-|R_0|$.  Since $K_0$ is not a coboundary, Observation~\ref{obs-nonsep}
implies that $\overline{K_0}$ contains a non-separating cycle.
Any non-separating cycle is non-contractible, and thus $|K_0|\ge \ew(G^\star)$.
Analogously, $|K_0+R_0|\ge \ew(G^\star)$.
Combining these relations and using Observation~\ref{obs-transv},
we conclude that
\begin{align*}
w(\PP_{G,f,Q})&\ge \frac{1}{2}\mu_f(R)\ge \frac{1}{2}\mu_f(R_0)\\
&=\frac{1}{2}(2\nu(R_0)+|R_0|-\trans{f}(R_0))\\
&=\frac{1}{2}(|K_0|+|K_0+R_0|-\trans{f}(R_0))\\
&=\frac{1}{2}\left(|K_0|+|K_0+R_0|-\sum_{v\in V(G)} \alpha_v(\partial_1 f)[v]\right)\\
&\ge \frac{1}{2}(|K_0|+|K_0+R_0|-|\partial_1 f|)\ge \ew(G^\star)-|\partial_1 f|/2.
\end{align*}
\end{proof}

We can now prove our result on 3-colorability of graphs with no odd-length faces.
\begin{proof}[Proof of Theorem~\ref{thm-genq}]
Suppose that $H$ has edgewidth at least $3\overline{\mu}_g$, and let $G$ be the dual of $H$.
Since every face of $H$ has even length, the dual graph $G$ is Eulerian, and thus it contains a nowhere-zero flow $f_0$
with $\partial_1 f_0=0$, obtained by sending a unit of flow along an Eulerian tour in $G$.  

Let $Q$ be a basis of $H^\star(G)$ and recall that $|Q|=g$.  For each $K\in Q$, let $r(K)=2\trans{f_0}(K)\bmod 3$.
Since $|\partial_1 f_0|=0$, Corollary~\ref{cor-width} implies that $w(\PP_{G,f_0,Q})\ge \ew(H)\ge 3\mu_g$,
and thus the polytope $\PP=(\PP_{G,f_0,Q}-r)/3$ has width at least $\ew(H)/3\ge \overline{\mu}_g$.
By Lemma~\ref{lemma-wide}, $\PP_{G,f_0,Q}$ is a translation of the polytope $\PP_{G,0,Q}$.
Observe that $\PP_{G,0,Q}=-\PP_{G,0,Q}$, and thus the polytopes $\PP_{G,0,Q}$, $\PP_{G,f_0,Q}$, and $\PP$
are centrally symmetric.  By Theorem~\ref{thm-narrow-censym}, this implies that $\PP$ contains an integer point.

By Theorem~\ref{thm-extcirc} and Observation~\ref{obs-polymod}, we conclude that
there exists an $f_0$-circulation $c$ such that $\trans{c}(K)\equiv r(K)\pmod 3$ for each $K\in Q$.
By Corollary~\ref{cor-circ}, there exists a nowhere-zero flow $f$ in $G$ with boundary $0$
such that $\trans{c}(K)\equiv 0\pmod 3$ for each $K\in Q$.
By Lemma~\ref{lemma-tutte}, this implies that $H$ is $3$-colorable.
\end{proof}

\section{Local 3-colorability}\label{sec-loc3col}

Finally, let us study local 3-colorability in graphs of large edgewidth.
First, let us note the following observation, which follows by the max-flow min-cut duality
and the correspondence between edge cuts and (edge-disjoint unions of) cycles in the dual graph.
Let $\Sigma$ be either an orientable surface other than the sphere, or the plane.
For a graph $H$ drawn in $\Sigma$, a $0$-boundary $d$ in the dual $G$ of $H$,
and a contractible cycle $K$ in $H$ bounding a disk $\Delta\subseteq\Sigma$,
let
$$\trans{d}(K)=\sum_{v\in V(G)\text{ drawn in $\Delta$}} d[v].$$
Let us remark that if $f$ is a flow in $G$ with $\partial_1 f=d$ and we view $K$ as a cocycle in $G$ in the
natural way, then $\trans{f}(K)=\trans{d}(K)$.
\begin{observation}\label{obs-cyccut}
Let $H$ be a graph with a 2-cell drawing either in an orientable surface other than the sphere or in the plane.
Let $G$ be the dual of $H$ and let $d$ be a $0$-boundary in $G$.  If $H$ is not drawn in the plane, suppose furthermore that
$\ew(H)\ge |d|/2$.  There exists a flow $f$ in $G$ with $\partial_1 f=d$
if and only if every contractible cycle $K$ in $H$ satisfies $|\trans{d}(K)|\le |K|$.
\end{observation}
Let $H$ be a graph drawn in an orientable surface of non-zero genus, and let $D(H)$ be the set of faces
of $H$ of length other than four.  For a non-empty set $L\subseteq D(H)$, we say a contractible cycle $K$ in $H$
\emph{surrounds} $L$ if $L$ is exactly the set of faces of $H$ of length other than four
drawn in the disk $\Delta\subseteq\Sigma$ bounded by $K$.  Let
$$b(L)=\min\Bigl(\sum_{x\in L} b(|x|),\sum_{x\in D(H)\setminus L} b(|x|)\Bigr),$$
and observe that
\begin{equation}\label{eq-lub}
b(L)\le \frac{b^\star(H)-1}{2}.
\end{equation}
We say that $L$ is \emph{dangerous} if $L$ is surrounded by a cycle in $H$ of length less than $b(L)$;
in this case, we let $K_L$ be a shortest cycle surrounding $L$, chosen arbitrarily,
and let $H_L$ be the subgraph of $H$ drawn in the closed disk bounded by $K_L$.
A single-element subset $L=\{x\}\subseteq D(H)$ is \emph{semi-dangerous} if no set containing $x$
is dangerous; in that case, we let $K_L$ be the facial walk of the face corresponding to $x$ and $H_L=K_L$.
A \emph{local 3-colorability witness} is the subgraph
$$U=\bigcup_{L\subseteq D(H):\text{ $L$ dangerous or semi-dangerous}} H_L$$
of $H$ (the naming is motivated by Lemma~\ref{lemma-planflow} below).
We need to show that local 3-colorability witnesses in graphs of sufficiently large edgewidth
are flat.  

\begin{lemma}\label{lemma-witflat}
Let $H$ be a simple graph with a 2-cell drawing in an orientable surface of non-zero genus and let $U\subseteq H$ be
a local 3-colorability witness.  If $\ew(H)\ge b^\star(H)$, then $U$ is flat.
\end{lemma}
\begin{proof}
We say a closed walk $W$ in $U$ is \emph{covered} by a sequence $L_1$, \ldots, $L_m$
of (not necessarily pairwise distinct) dangerous or semi-dangerous sets if $W$ is a concatenation of walks
$W_1$, \ldots, $W_m$, where $W_i$ is a walk in $H_{L_i}$
for $i\in\{1,\ldots,m\}$.
Suppose for a contradiction that $U$ contains a non-contractible walk $W$, and let us choose one covered by
a sequence $L_1$, \ldots, $L_m$ of dangerous or semi-dangerous sets with $m$ smallest possible.

Let $C$ be a cycle with vertices $1$, \ldots, $m$ in order.  Consider any indices $i<j$ non-adjacent in $C$.
We claim that $L_i\cap L_j=\emptyset$.  Indeed, if $x\in L_i\cap L_j$, then the facial walk of $x$
is contained in both $H_{L_i}$ and $H_{L_j}$, and since these subgraphs are connected,
there exists a path $W_0$ from a vertex $z_i$ of $W_i$ to a vertex $z_j$ of $W_j$ such that $W_0$ is a concatenation of
a path in $H_{L_i}$ and a path in $H_{L_j}$.  For $a\in \{i,j\}$, express $W_a$ as the concatenation of walks $W'_a$
and $W''_a$ ending and starting in $z_a$, respectively.  Let $W'$ be the concatenation of $W''_j$, $W_{j+1}$, \ldots, $W_m$, $W_1$,
\ldots, $W_{i-1}$, $W'_i$, and $W_0$.  Let $W''$ be the concatenation of $W''_i$, $W_{i+1}$, \ldots, $W_{j-1}$, $W'_j$, and the
reverse of $W_0$, see Figure \ref{fig-LemWitFlatCon}.  Since the closed walk $W$ is non-contractible, by symmetry between $W'$ and $W''$ we can assume that
the closed walk $W''$ is non-contractible.
However, $W''$ is covered by $L_i$, $L_{i+1}$, \ldots, $L_j$, contradicting the choice of $W$.

\begin{figure}
\centering{
\begin{tikzpicture}[scale=.9]
\draw (0,1.5) ellipse (1cm and 2cm);
\draw (0,-1.5) ellipse (1cm and 2cm);
\draw[rotate=10] (-2, 3) ellipse (2cm and 1cm);
\draw[rotate=170] (2, 3) ellipse (2cm and 1cm);
\draw (-4, 0) ellipse (1cm and 2cm);
\draw[rotate=-10] (2, 3) ellipse (2cm and 1cm);
\draw[rotate=-170] (-2, 3) ellipse (2cm and 1cm);
\draw (4, 0) ellipse (1cm and 2cm);
\draw (-.1, 0) -- (-.1, .25) -- (.2, .3) -- (.4, .1) -- (.2, -.3) -- (-.1, 0);
\draw[thick, blue] plot [smooth] coordinates {(0, 2.3) (-1.1, 2.7) (-2, 3) (-3.5, 2.5) (-4.2, .8) (-3.5, -.4) (-3.8, -2) (-3, -2.5) (-2, -3) (-.7, -2.4) (0, -2.1)};
\draw[thick, violet] plot [smooth] coordinates {(0, 2.3) (.2, 1.3) (-.4, 0)  (0, -2.1)};
\draw[thick, red] plot [smooth] coordinates {(0, -2.1) (1.1, -2.2) (3.1, -3) (3.6, -1) (4.1, 2) (2.8,3.3) (1.7,  2.8) (0,2.3)};
\fill[black] (0, 2.3) circle(1.5pt);
\fill[black] (0, -2.1) circle(1.5pt);
\node at (0, 2.5) {$z_i$};
\node at (0, -2.4) {$z_j$};
\node at (-3.3, 3.8) {$K_{L_{i+1}}$};
\node at (0, 3.8) {$K_{L_i}$};
\node at (3.3, 3.8) {$K_{L_{i-1}}$};
\node at (-3.3, -3.8) {$K_{L_{j-1}}$};
\node at (0, -3.8) {$K_{L_j}$};
\node at (3.3, -3.9) {$K_{L_{j+1}}$};
\node[blue] at (-2.4, 2.5) {$W''$};
\node[red] at (2.6, 2.8) {$W'$};
\node[violet] at (.5, 1) {$W_0$};
\node at (.13, .1) {$x$};
\end{tikzpicture}}
\caption{The construction of $W'$ and $W''$ in Lemma \ref{lemma-witflat}.}
\label{fig-LemWitFlatCon}
\end{figure}
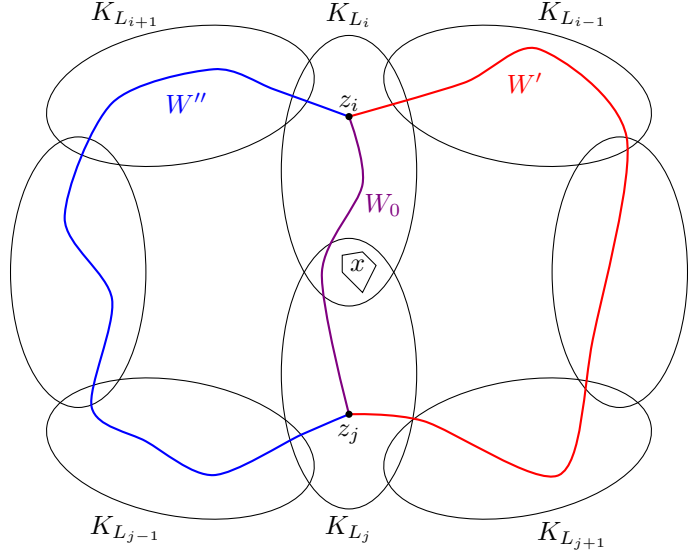

Since $W$ is non-contractible, we have $W\not\subseteq H_{L_1}$, and thus $m>1$.  Let $Y$ be the closed
walk equal to concatenation of walks $Y_1$, \ldots, $Y_m$ obtained as follows: We start by setting $Y=W$ and $Y_i=W_i$
for each $i$.  Then, for $i=1,\ldots,m$:
\begin{itemize}
\item Replace $Y_i$ by the longest subwalk of $Y$ contained in $H_{L_i}$, and shorten $Y_{i-1}$ and $Y_{i+1}$
correspondingly (note that they stay non-empty by the minimality of $m$ from the choice of $W$).
Observe that both ends $x_i$ and $y_i$ of $Y_i$ are contained in $K_{L_i}$.
\item Replace $Y_i$ by the shorter of the two walks between $x_i$ and $y_i$ in $K_{L_i}$, see Figure \ref{fig-LemWitFlatYs}.  Note that since $Y_i$
is contained in the disk bounded by $K_{L_i}$, $Y$ remains homotopically equivalent to $W$.
\end{itemize}
At the end, $Y$ is a non-contractible closed walk in $H$, and for each $i$, the length of $Y_i$ is at most half the length of $K_i$.
Since $Y$ is non-contractible, the bound on the edgewidth of $H$ gives $|Y|\ge \ew(H)\ge b^\star(H)$.

\begin{figure}
\centering{
\begin{tikzpicture}[scale = .9]
\draw (2, -.05) arc(210:-60:2cm and 1cm);
\draw (-1, -.3) arc(230:-31:2cm and 1cm);
\draw(-4.9, 0) arc(180:-17:2cm and 1cm);
\draw[red,thick] (2, -.05) arc(210:300:2cm and 1cm);
\draw[red,thick] (-1, -.3) arc(230:329:2cm and 1cm);
\draw[red,thick] (-4.9, 0) arc(180:343:2cm and 1cm);
\draw[thick, blue] plot [smooth] coordinates {(-4.9,0) (-4,.5) (-3.4, .6) (-3.2, .1) (-2, -.4) (-1.3, -.1) (0, .4) (.6, .3) (1.9, 0) (2.8, .7) (3.9, .5) (4.4, -.1) (5.05, -.3)};
\node at (-3, -.3) {$W_{i+1}$};
\node at (.6, .6) {$W_{i}$};
\node at (4.4, .6) {$W_{i-1}$};
\node at (-2.9, -1.3) {$Y_{i+1}$};
\node at (.5, -.8) {$Y_{i}$};
\node at (3.5, -.9) {$Y_{i-1}$};
\node at (-3, 1.3) {$K_{L_{i+1}}$};
\node at (.2, 1.8) {$K_{L_i}$};
\node at (3.7, 1.8) {$K_{L_{i-1}}$};
\end{tikzpicture}}
\caption{Each walk $W_i$ is replaced by a walk $Y_i$ along the cycle $K_{L_i}$.}
\label{fig-LemWitFlatYs}
\end{figure}
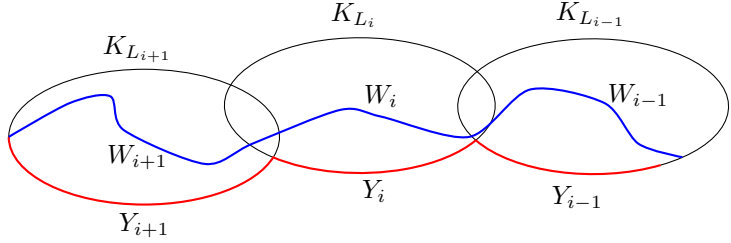

Let $S\subseteq D(H)$ consist of the faces $x\in D(H)$ such that $\{x\}$ is semi-dangerous, and let
$$\sigma=\sum_{x\in S} b(|x|).$$
For any dangerous set $L$, we have $L\cap S=\emptyset$, and thus
$$b(L)\le \sum_{x\in L} b(|x|)\le b^\star(H)-\sigma-1.$$
Together with (\ref{eq-lub}), this gives that for each $i$, 
$$|Y_i|\le \Bigl\lfloor \frac{|K_i|}{2}\Bigr\rfloor\le
\begin{cases}
\frac{b(L_i)}{2}\le \min\Bigl(\frac{b^\star(H)-1}{4},\frac{b^\star(H)-1-\sigma}{2}\Bigr)&\text{ if $L_i$ is dangerous}\\
\lfloor |x|/2\rfloor \le b(|x|)&\text{ if $L_i=\{x\}$ is semi-dangerous.}
\end{cases}$$
Note that for each $x\in D(H)$, if $\{x\}$ is semi-dangerous, then
there exists at most one index $i$ such that $x\in L_i$, and we have $L_i=\{x\}$ for this index
(if $x\in L_j$ for $j\neq i$, then $j$ would be adjacent to $i$ in $C$ and $L_j=\{x\}=L_i$,
contradicting the minimality of $m$).  Consequently,
$$\sum_{i:L_i\text{ semi-dangerous}} |Y_i|\le \sum_{i:L_i=\{x_i\}\text{ semi-dangerous}} b(|x_i|)\le\sigma.$$
Distinguishing two cases depending on whether $L_i$ is dangerous for all $i\in\{1,\ldots,m\}$ or not, we conclude that if $m\le 3$, then
$$|Y|=\sum_{i=1}^m |Y_i|\le \max\left(m\cdot \frac{b^\star(H)-1}{4}, (m-1)\cdot \frac{b^\star(H)-1-\sigma}{2}+\sigma\right)<b^\star(H).$$
This is a contradiction, and thus $m\ge 4$.  Since $L_i\cap L_j=\emptyset$ whenever $i$ and $j$ are non-adjacent in $C$,
for each $x\in D(H)$ such that $\{x\}$ is not semi-dangerous, there exist at most two indices $i$ such that $x\in L_i$.  
Hence,
\begin{align*}
|Y|&=\sum_{i=1}^m |Y_i|\le \sigma+\sum_{i:L_i\text{ dangerous}} |Y_i|\le \sigma+\frac{1}{2}\sum_{i:L_i\text{ dangerous}} b(L_i)\\
&\le \sigma+\frac{1}{2}\sum_{i:L_i\text{ dangerous}} \sum_{x\in L_i} b(|x|)
\le \sigma+\frac{1}{2}\sum_{x\in D(H)\setminus S} 2\cdot b(|x|)=b^\star(H)-1.
\end{align*}
This is again a contradiction.
\end{proof}

Let $U$ be a flat subgraph of a graph $H$ drawn in a surface.  A subgraph $U_0\subseteq H$ is a \emph{consolidation}
of $U$ if it is obtained from $U$ by repeatedly adding paths (in $H$) between distinct connected components until $U_0$ is connected.
Note that every cycle in $U_0$ is also a cycle in $U$, and thus $U_0$ is also flat.

Consider a connected flat subgraph $U$ of $H$, let $U'$ be a planar quadrangulation extension of $U$, and let $\theta$ be
the homeomorphism mapping $U$ to a subgraph of $U'$ from the definition of planar quadrangulation extension.  Let $W$
be the closed walk bounding the outer face of $U$, and let $W'=\theta(W)$ be the corresponding closed walk in $U'$.
A walk $P$ in $U'$ is a \emph{chord} if its ends are in $W'$, all other vertices and edges of $P$ are drawn in the outer face of $\theta(U)$,
and all vertices of $P$ except possibly for its ends are pairwise distinct
(i.e. $P$ is either a path with both ends in $W'$, or a cycle intersecting $W'$ in exactly one vertex).  The plane graph $\theta(U)+P$ has exactly
one internal face that is not a face of $\theta(U)$, and the facial walk of this face consists of $P$ and a subwalk $B_P$
of $W$ or the reverse of $W$; we say that $B_P$ is the \emph{base} of the chord.  We say that the
planar quadrangulation extension $U'$ of $U$ is \emph{generic} if each chord is at least as long as its base.

Observe that if it is possible to quadrangulate the outer face of a connected flat subgraph to obtain a planar quadrangulation extension,
it is also possible to do it so that the resulting planar quadrangulation extension is generic.
\begin{observation}\label{obs-ex}
Let $H$ be a graph with a 2-cell drawing in an orientable surface of non-zero genus and let $U$ be a connected flat subgraph of $H$.
Then the following claims are equivalent:
\begin{itemize}
\item the outer face of $U$ has even length,
\item $U$ has a generic planar quadrangulation extension,
\item $U$ has a planar quadrangulation extension.
\end{itemize}
\end{observation}

Let us now give the key property of a local 3-colorability witness.

\begin{lemma}\label{lemma-planflow}
Let $H$ be a simple graph with a 2-cell drawing in an orientable surface of non-zero genus with edgewidth at least $b^\star(H)/2$,
let $G$ be the dual of $H$, and let $U\subseteq H$ be a flat subgraph capturing non-$4$-faces of $H$.
\begin{itemize}
\item If $H$ is $3$-colorable, then every generic planar quadrangulation extension $U_1$ of a consolidation $U_0$ of $U$ is $3$-colorable.
\item If $U$ is a local 3-colorability witness and a planar quadrangulation extension of $U$ is $3$-colorable, then
there exists a nowhere-zero flow $f$ in $G$ with $3|\partial_1 f$.
\end{itemize}
\end{lemma}
\begin{proof}
Recall that $D(H)$ is the set of faces of $H$ of length other than four, and let $I$ be the set of corresponding vertices of $G$.
Since $U$ captures non-$4$-faces of $H$, each face in $D(H)$ is also an internal face of $U$.

For the first claim, we translate a 3-coloring of $H$ using
Lemma~\ref{lemma-tutte} to a nowhere-zero flow $f_0$ in $G$, and we copy its
boundary $d_0$ to a boundary $d_1$ in $U_1^\star$.  We then use
Observation~\ref{obs-cyccut} to show that there exists a flow $f_1$ in $U_1^\star$ with boundary $d_1$,
extend $f_1$ to a nowhere-zero flow $f_1'$, and finish by converting $f_1'$ to a 3-coloring
of $U_1$ using Lemma~\ref{lemma-tutte}.  Let us now state this argument precisely.

Let $\theta$ be a homeomorphism mapping $U_0$ to a subgraph of $U_1$ from the definition of planar quadrangulation
extension.  Slightly abusing the notation, let us also view $\theta$ as mapping vertices of $I$ to the corresponding vertices of the plane dual
$U^\star_1$ of $U_1$.
Since $H$ is 3-colorable, Lemma~\ref{lemma-tutte} implies that there exists a nowhere-zero flow $f_0$
in $G$ with boundary $d_0=\partial_1 f_0$ divisible by $3$.  Let $d_1$ be the $0$-boundary in $U^\star_1$ defined as
$$d_1=\sum_{v\in I} d_0[v]\cdot \theta(v).$$
Consider any cycle $K$ in $U_1$, and let $M$ be the set of vertices of $\theta(I)$ drawn in the disk in the plane bounded by $K$.
If $M=\theta(I)$ or $M=\emptyset$, then since $d_1$ is a $0$-boundary, we have $\trans{d_1}(K)=0\le |K|$.  Otherwise,
observe that $K$ must intersect $\theta(U_0)$. Let $K'$ be the closed walk obtained from $K$ by replacing all chords by their
bases; since $U_1$ is generic, we have $|K'|\le |K|$.  Note that $K'\subseteq \theta(U_0)$, and that
by the definition of the base of a chord, the set of vertices of $\theta(I)$ drawn in the region bounded by
$K'$ is the same as for $K$, i.e., $M$.  Observe that $\theta^{-1}(M)$ is exactly the set
of vertices of $I$ drawn in the part of $\Sigma$ bounded by $\theta^{-1}(K')$.  Therefore,
viewing $\theta^{-1}(K')$ as a cocycle in $G$, we have
$$\trans{d_1}(K)=\sum_{v\in M} d_1[v]=\sum_{u\in\theta^{-1}(M)} d_0[u]=\trans{f_0}(\theta^{-1}(K'))\le |K'|\le |K|.$$
Since this holds for every cycle $K$, Observation~\ref{obs-cyccut} implies that there exists
a flow $f_1$ in $U^\star_1$ with $\partial_1 f_1=d_1$.  Since $d_0$ is parity-compliant by Lemma~\ref{lemma-flow1},
$d_1$ is also parity-compliant, and thus Lemma~\ref{lemma-flow1} implies that there exists
a nowhere-zero flow $f'_1$ in $U^\star_1$ with $\partial_1 f'_1=d_1$.  Since $3|d_0$, we also have $3|d_1$,
and thus Lemma~\ref{lemma-tutte} implies that $U_1$ is $3$-colorable.

\medskip

Suppose now that $U$ is a local 3-colorability witness and a planar
quadrangulation extension $U_2$ of $U$ is $3$-colorable. To prove the second
claim, we use Lemma ~\ref{lemma-tutte} to obtain a nowhere-zero flow $f'$ in
$U_2^\star$ with boundary $d'$ divisible by $3$. We then translate $d'$ into a $0$-boundary $d$
in $G$ and show by contradiction through Observation~\ref{obs-ex} that there
exists a flow $f$ in $G$ with $\partial_1 f=d$.

More precisely, let $\theta$ be a homeomorphism mapping $U$ to a subgraph of $U_2$ from the definition of planar quadrangulation extension.
Since $U_2$ is 3-colorable, by Lemma~\ref{lemma-tutte} there exists a nowhere-zero flow $f'$ in $U^\star_2$ with boundary $d'=\partial_1 f'$
divisible by $3$.  Since $U$ is a local 3-colorability witness, all internal faces of $U$ are also
faces of $H$ and all other faces of $H$ have length four.  Let
$$d=\sum_{v\in I} d'[\theta(v)]\cdot v;$$
then $d$ is a parity-compliant $0$-boundary in $G$ such that $3|d$, and thus
by Lemma~\ref{lemma-flow1}, it suffices to show that $G$ contains a flow with $\partial_1 f=d$.  Suppose for a contradiction that this is not the
case.  Note that $|d|/2<b^\star(H)/2\le \ew(H)$, and thus by Observation~\ref{obs-cyccut}, there exists a cycle $K$ in $H$ bounding
a disk $\Delta\subseteq\Sigma$ such that
\begin{equation}\label{eq-shk}
\trans{d}(K)>|K|.
\end{equation}
For any face $x$ of $H$ of length other than four, if $v\in I$ is the corresponding vertex of $G$, then
$\theta(v)$ is a vertex of $U^\star_2$ of degree $\deg v=|x|$,
and by the existence of the flow $f'$ and the fact that $3|d'$ and $d'$ is parity-compliant,
we have $|d[v]|=|d'[\theta(v)]|\le b(\deg v)=b(|x|)$.

Let $L$ be the set of faces of $H$ drawn in $\Delta$ of length other than four, and let $L^\star$ be the set of corresponding vertices of $G$.
Note that
$$|\trans{d}(K)|=\Bigl|\sum_{v\in L^\star} d[v]\Bigr|\le \sum_{x\in L} b(|x|),$$
and since $d$ is a $0$-boundary, we also have
$$|\trans{d}(K)|=\Bigl|\sum_{v\in I\setminus L^\star} d[v]\Bigr|\le \sum_{x\in D(H)\setminus L} b(|x|).$$
Consequently, $|\trans{d}(K)|\le b(L)$, and (\ref{eq-shk}) implies that
$|K|<b(L)$. That is, $L$ is dangerous.  Since $U$ is a local 3-colorability witness, there
exists a cycle $K'$ in $U$ of length at most $|K|$ surrounding $L$.  However, then observe that $\theta(L)$ is
exactly the set of faces of $U_2$ of length other than four drawn in the cycle $\theta(K')\subseteq U_2$ and
$$\trans{d'}(\theta(K'))=\trans{d}(K)>|K|\ge |\theta(K')|.$$
By Observation~\ref{obs-cyccut}, this contradicts the existence of the flow $f'$.
\end{proof}

We are ready to prove the characterization of 3-colorability for graphs of large edgewidth.
\begin{proof}[Proof of Theorem~\ref{thm-ew}]
If $H$ is 3-colorable, then by Observation~\ref{obs-ex}, a consolidation of any flat subgraph capturing non-$4$-faces
has a generic planar quadrangulation extension $U_1$, and by Lemma~\ref{lemma-planflow}, $U_1$ is 3-colorable.
It follows that $H$ is locally 3-colorable.

Conversely, suppose that $H$ is locally 3-colorable and its edgewidth is at least $b^\star(H)-1+3\mu_g$.
Let $U$ be a local 3-colorability witness; by Lemma~\ref{lemma-witflat}, $U$ is flat.  
Furthermore, since each face of $H$ of length other than
four either forms a semi-dangerous set or is contained in a dangerous set,
$U$ captures non-$4$-faces of $H$
Let $G$ be the dual of $H$.  Since $H$ is locally 3-colorable, there exists a $3$-colorable planar quadrangulation extension of $U$,
and by Lemma~\ref{lemma-planflow}, there exists a nowhere-zero flow $f_0$ in $G$ with $3|\partial_1 f_0$.
Let $Q$ be a basis of $H^\star(G)$, and recall that $|Q|=g$.  For each $K\in Q$, let $r(K)=2\trans{f_0}(K)\bmod 3$.
By Corollary~\ref{cor-width}, we have $w(\PP_{G,f,Q})\ge \ew(H)-|\partial_1 f_0|/2\ge \ew(H)-(b^\star(H)-1)\ge 3\mu_g$.
Consequently, the polytope $(\PP_{G,f_0,Q}-r)/3$ has width at least $\mu_g$, and thus by Theorem~\ref{thm-narrow},
it contains an integer point.  By Theorem~\ref{thm-extcirc} and Observation~\ref{obs-polymod}, there
exists an $f_0$-circulation $c$ such that $\trans{c}(K)=r(K)\pmod 3$ for each $K\in Q$.
By Lemma~\ref{lemma-tutte} and Corollary~\ref{cor-circ}, this implies that $H$ is $3$-colorable.
\end{proof}

The improved bound for the torus is proved analogously.
\begin{proof}[Proof of Corollary~\ref{cor-torus}]
We proceed exactly as in the second part of the proof of Theorem~\ref{thm-ew}.  Since the polytope $(\PP_{G,f_0,Q}-r)/3$ is $\tfrac{1}{3}$-integral
and $g=2$, by Lemma~\ref{lemma-narrow} it suffices to show that its width is at least two, which is implied by the assumption
$\ew(H)\ge b^\star(H)+5$.
\end{proof}

\section{Concluding remarks}\label{sec-final}

Let us note that no superlinear lower bound on the constant $\mu_d$ from Theorem~\ref{thm-narrow} is known, and indeed it has been conjectured that $\mu_d=O(d\log d)$ or even $\mu_d=O(d)$.
Any improvement over the current $O(d^{4/3})$ bound directly translates into improved bounds for Theorem~\ref{thm-ew}.

It is natural to ask whether our techniques also apply to non-orientable surfaces.  While this seems
to be the case to some extent, there are additional challenges coming from the fact that the homology group of a non-orientable
surface of Euler genus $g$ is $\mathbb{Z}^{g-1}\times \mathbb{Z}_2$, leading to parity considerations that do not arise in
the orientable case.  Dvořák, Moore and Sereni (private communication) worked out these issues for the projective plane and Klein bottle,
obtaining simple efficient algorithms for these surfaces.

In Lemma~\ref{lemma-wide}, we prove that the polytope $\PP_{G,f,Q}$ is a translation of the polytope $\PP_{G,\partial_1 f,Q}$,
and thus it is essentially independent of the exact choice of the initial nowhere-zero flow, subject to its fixed boundary. 
It is actually possible to eliminate the dependence on the boundary and make the connection to coloring even more direct.
Let $L$ be the set of vertices of $G$ of degree other than four.  
We define
$$\PP_{G,Q}=\left\{\begin{array}{l}
(a,d)\in \mathbb{R}^Q\times\mathbb{R}^L: \langle d,1\rangle=0,\\
\langle z,a\rangle + \langle z',d\rangle 
\le \min_{z''\in \mathbb{Z}^{V(G)\setminus L}}|\langle z,Q\rangle + \langle (z',z''),\partial^\star_2\rangle|\\
\text{\hspace{26mm}for every $(z,z')\in \mathbb{Z}^Q\times \mathbb{Z}^L$}
\end{array}\right\}.$$
For $K\in Q$, let $a_0(K)=3$ if $|K|$ is odd and $a_0(K)=0$ if $|K|$ is even.
For $v\in L$, let $d_0(v)=3$ if $\deg v$ is odd and $d_0(v)=0$ if $\deg v$ is even.
Applying the theory developed in this paper, it is easy to see that $H=G^\star$ is $3$-colorable if and only if $(\PP_{G,Q}-(a_0,d_0))/6$ contains an integer point.
Hence, using integer programming in bounded dimension, we actually obtain a polynomial-time algorithm for 3-coloring
graphs drawn in a fixed orientable surface as long as the number of non-4-faces is bounded (without any restrictions on their length).
However, in the case $H$ is a near-quadrangulation, it is likely faster to go through at most $q^\star(H)$ plausible boundaries
and apply the integer programming to the polytopes $\PP_{G,f,Q}$ rather than to the higher-dimensional polytope $\PP_{G,Q}$.

Moreover, the polytope $\PP_{G,Q}$ is centrally symmetric; one could hope to use this fact in connection with Theorem~\ref{thm-narrow-censym}
to improve the bound in Theorem~\ref{thm-ew}.  However, the issue here is that the width of $\PP_{G,Q}$ cannot be lower bounded in
terms of the edgewidth of $G$, since it can be narrow in a direction $(0,z')$ for some $z'\neq 0$---in this case, the constraints are
only affected by coboundaries.

\bibliographystyle{plainurl}
\bibliography{references}

\end{document}